\documentclass[a4paper,12pt]{article}

\usepackage{amsmath}
\usepackage{ntheorem}
\usepackage{amssymb}
\usepackage{latexsym}
\usepackage{rsfs}
\usepackage[all]{xy}\UseComputerModernTips
\usepackage{graphicx}

\theoremstyle{plain}
\newtheorem{prp}{Proposition}[section]
\newtheorem{thm}[prp]{Theorem}
\newtheorem{lem}[prp]{Lemma}
\newtheorem{cor}[prp]{Corollary}

\theoremstyle{plain}
\theorembodyfont{\normalfont\rm}
\newtheorem{dfn}[prp]{Definition}
\newtheorem{rem}[prp]{Remark}
\newtheorem{exm}[prp]{Example}

\theoremstyle{nonumberplain}
\theoremheaderfont{\normalfont\itshape}
\theoremseparator{.}
\theorembodyfont{\normalfont}
\newtheorem{proof}{Proof}

\newcommand{\qed}{\hfill $\Box$}

\makeatletter

\@addtoreset{equation}{section}
\makeatother

\newcommand{\Real}{{\mathbb R}}
\newcommand{\Comp}{{\mathbb C}}
\newcommand{\Z}{{\mathbb Z}}
\renewcommand{\P}{{\mathbb P}}

\newcommand{\D}{{\rm D}}
\newcommand{\Db}{{\bf D}^{b}}
\newcommand{\Dc}{{\bf D}_{\Real -c}^{b}}

\newcommand{\tLb}{{\cal E}}
\newcommand{\Lb}{{\cal F}}
\newcommand{\sLb}{{\cal G}}
\newcommand{\Vb}{{\cal G}}
\newcommand{\B}{{\cal B}}

\renewcommand{\H}{{\cal H}}
\renewcommand{\L}{{\mathscr L}}

\renewcommand{\1}{{\bf 1}}
\newcommand{\CF}{{\rm CF}}
\renewcommand{\S}{{\mathfrak S}}

\newcommand{\reg}{{\rm reg}}
\newcommand{\tr}{{\rm tr}}
\newcommand{\id}{{\rm id}}
\newcommand{\supp}{{\rm supp}}
\renewcommand{\SS}{{\rm SS}}
\newcommand{\Ev}{{\rm Ev}}
\newcommand{\sgn}{{\rm sgn}}
\newcommand{\RG}{R\varGamma}

\newcommand{\Hom}{{\cal H}om}
\newcommand{\e}{\varepsilon}

\renewcommand{\t}[1]{{^{t}{#1}^{\prime}}}
\newcommand{\tl}[1]{\widetilde{#1}}
\newcommand{\uotimes}[1]{\otimes_{#1}}
\newcommand{\utimes}[1]{\times_{#1}}
\renewcommand{\d}{{\rm dim}}

\newcommand{\simto}{\overset{\sim}{\longrightarrow}}
\newcommand{\simot}{\overset{\sim}{\longleftarrow}}
\newcommand{\dsum}{\displaystyle \sum}
\newcommand{\dint}{\displaystyle \int}
\renewcommand{\injlim}[1]{\underset{#1}{\varinjlim}}
\renewcommand{\(}{\left(}
\renewcommand{\)}{\right)}
\newcommand{\inun}{\text{\rotatebox[origin=c]{90}{$\in$}}}
\newcommand{\longhookrightarrow}{\DOTSB\lhook\joinrel\longrightarrow}
\newcommand{\longtwoheadrightarrow}{\relbar\joinrel\twoheadrightarrow}

\title{Microlocal study of Lefschetz fixed point formulas for higher-dimensional fixed point sets}

\author{Yutaka Matsui\footnote{Department of Mathematics, Kinki University, 3-4-1, Kowakae, Higashi-Osaka, Osaka, 577-8502, Japan, E-mail: matsui@math.kindai.ac.jp} \and Kiyoshi Takeuchi\footnote{Institute of Mathematics, University  of Tsukuba, 1-1-1, Tennodai, Tsukuba, Ibaraki, 305-8571, Japan, E-mail: takemicro@nifty.com} }
\date{ }

\begin{document}

\maketitle

\begin{abstract}
We introduce new Lagrangian cycles which encode local contributions of Lefschetz numbers of constructible sheaves into geometric objects. We study their functorial properties and apply them to Lefschetz fixed point formulas with higher-dimensional fixed point sets.
\end{abstract}

\section{Introduction}
The aim of this paper is to study Lefschetz fixed point formulas for morphisms $\phi \colon X \longrightarrow X$ of real analytic manifolds $X$ whose fixed point set $M= \{ x \in X \ |\ \phi(x)=x\} \subset X$ is higher-dimensional (since we mainly consider the case where the fixed point set is a smooth submanifold of $X$, we use the symbol $M$ to express it). Since the beginning of the theory, it is well-known that when $X$ is compact the global Lefschetz number of $\phi$
\begin{equation}
\tr(\phi):=\dsum_{j \in \Z} (-1)^j \tr\{H^j(X;\Comp_X)\overset{\phi^*}{\longrightarrow} H^j(X;\Comp_X) \} \in \Comp
\end{equation}
is expressed as the integral of a local cohomology class $C(\phi) \in H_M^n(X;or_X)$ supported by $M$, where we set $\d X=n$ and $or_X$ is the orientation sheaf of $X$. See e.g. Dold \cite{Dold}, \cite{Dold-2} etc. for the detail of this subject. Let $M= \bigsqcup_{i \in I}M_i$ be the decomposition of $M$ into connected components and 
\begin{gather}
H_M^n(X;or_X) =\bigoplus_{i \in I} H_{M_i}^n (X;or_X),\\
C(\phi) =\bigoplus_{i \in I} C(\phi )_{M_i}
\end{gather}
the associated direct sum decompositions. Then the integral of the local cohomology class $C(\phi)_{M_i} \in H_{M_i}^n (X;or_X)$ associated with a fixed point component $M_i$ is called the local contribution from $M_i$. In other words, the global Lefschetz number of $\phi$ is equal to the sum of local contributions from $M_i$'s. But if the fixed point component $M_i$ is higher-dimensional, it is in general very difficult to compute the local contribution by the following dimensional reason. Let $M_i$ be a fixed point component of $\phi$ whose codimension in $X$ is $d>0$. Then the local cohomology group $H_{M_i}^n (X;or_X)$ is isomorphic to the $0$-dimensional Borel-Moore homology group $H^{\text{B-M}}_0(M_i;\Comp)$ of $M_i$ by the Alexander duality, and the class $C(\phi)$ in it is very hard to deal with. Recall that top-dimensional Borel-Moore homology cycles in $M_i$, i.e. elements in $H^{\text{B-M}}_{n-d}(M_i;\Comp)$, can be much more easily studied since they are realized as sections of a relative orientation sheaf on $M_i$.

In this paper, we overcome this difficulty for smooth fixed point components $M_i$ by introducing new Lagrangian cycles in the cotangent bundles $T^*M_i$. Since we also want to study Lefschetz fixed point formulas over singular varieties (and those for intersection cohomology groups), we consider the following very general setting. Let $X$, $\phi$ and $M= \bigsqcup_{i \in I}M_i$ be as before, and $F$ a bounded complex of sheaves of $\Comp_X$-modules whose cohomology sheaves are $\Real$-constructible in the sense of \cite{K-S}. Assume that we are given a morphism $\Phi \colon \phi^{-1}F \longrightarrow F$ in the derived category $\Dc(X)$. If the support $\supp(F)$ of $F$ is compact, we can define the global trace (Lefschetz number) $\tr (F,\Phi) \in \Comp$ of the pair $(F, \Phi)$ by
\begin{equation}
\tr (F,\Phi) :=\dsum_{j \in \Z} (-1)^j \tr \{H^j(X;F)\overset{\Phi}{\longrightarrow} H^j(X;F) \} \in \Comp,
\end{equation}
where the morphisms $H^j(X;F) \overset{\Phi}{\longrightarrow} H^j(X;F)$ are induced by
\begin{equation}
F \longrightarrow R\phi_*\phi^{-1}F \overset{\Phi}{\longrightarrow}R\phi_*F.
\end{equation}
Also in this very general setting, Kashiwara \cite{Kashiwara-2} introduced local contributions $c(F,\Phi)_{M_i} \in \Comp$ from the fixed point components $M_i$ and proved that
\begin{equation}
\tr (F,\Phi)=\dsum_{i \in I} c(F,\Phi)_{M_i}.
\end{equation}
Therefore the next important problem in the Lefschetz fixed point formula for constructible sheaves is to describe these local contributions $c(F,\Phi)_{M_i}$.

Now let us take a smooth fixed point component $M_i$. For the sake of simplicity, we shall denote it also by $M$. Then there exists a natural morphism
\begin{equation}
\phi^{\prime} \colon T_MX \longrightarrow T_MX
\end{equation}
induced by $\phi$, where $T_MX$ is the normal bundle of $M$ in $X$. For each point $x\in M$, we define a finite subset $\Ev_x$ of $\Comp$ by
\begin{equation}
\Ev_x:=\{ \text{the eignevalues of $\phi^{\prime}_x \colon (T_{M}X)_x \longrightarrow (T_{M}X)_x$} \} \subset \Comp.
\end{equation}
Assume the condition:
\begin{equation}\label{eq:1-9}
1 \notin \Ev_x \hspace{5mm} \text{for any $x \in \supp(F) \cap M$},
\end{equation}
which means that the graph $\Gamma_{\phi}=\{ (\phi(x), x) \ | \ x\in X\} \subset X \times X$ of $\phi$ intersects cleanly (see \cite[Definition 4.1.5]{K-S}) with the diagonal set $\Delta_X$ along $M \subset \Gamma_{\phi} \cap \Delta_X$ on $\supp(F) \cap M$. This condition naturally appears also in the study of Atiyah-Bott type (holomorphic) Lefschetz theorems by Toledo-Tong \cite{T-T-3}. Under the condition \eqref{eq:1-9}, in Section \ref{sec:4} we shall construct a new Lagrangian cycle $LC(F, \Phi)_M$ in the cotangent bundle $T^*M$. In this paper, we call this cycle $LC(F, \Phi)_M$ the Lefschetz cycle associated with the pair $(F, \Phi)$ and the fixed point component $M$. Note that in the more general setting of elliptic pairs a similar construction of microlocal Lefschetz classes was already given in the pioneering work \cite{Guillermou} by Guillermou. The difference from his construction is that we explicitly realize such microlocal characteristic classes as geometric objects in the cotangent bundle $T^*M$. Note also that if $\phi=\id_X$, $M=X$ and $\Phi =\id_F$, our Lefschetz cycle $LC(F, \Phi)_M$ coincides with the characteristic cycle $CC(F)$ of $F$ introduced by Kashiwara \cite{Kashiwara-1}. For the applications of characteristic cycles to projective duality, see \cite{Ernstrom}, \cite{M-T-2}, \cite{M-T-1} etc. By Lefschetz cycles, we can generalize almost all nice properties of characteristic cycles into more general situations.

First, in Section \ref{sec:5} we prove the following microlocal index theorem for the local contribution $c(F, \Phi)_M$ from $M$ (see Theorem \ref{thm:5-1}).

\begin{thm}
Assume that $\supp (F) \cap M$ is compact. Then for any continuous section $\sigma \colon M \longrightarrow T^*M$ of $T^*M$, we have
\begin{equation}
c(F,\Phi)_M =\sharp ([\sigma] \cap LC(F,\Phi)_M),
\end{equation}
where $\sharp ([\sigma] \cap LC(F,\Phi)_M)$ is the intersection number of the image of $\sigma$ and $LC(F,\Phi)_M$ in the cotangent bundle $T^*M$.
\end{thm}

Next in Section \ref{sec:6}, by using this microlocal index theorem, in many cases we give some useful formulas, similar to those for characteristic cycles, for the explicit description of our Lefschetz cycles. In the course of the proof of these results, we obtain in Section \ref{sec:3} and \ref{sec:6} some localization theorems which partially generalize the previous results by Goresky-MacPherson \cite{G-M-1}, Kashiwara-Schapira \cite{K-S} and Braden \cite{Braden} etc. In particular, to prove the localization theorem in the case where the set $\Ev_x$ may vary depending on $x \in \supp(F) \cap M$, we required some precise arguments on Lefschetz cycles in Section \ref{sec:6} (see also Remark \ref{rem:3-6}). As we shall see in Example \ref{exm:3-8} and \ref{exm:3-9}, these localization theorems also have some applications to the explicit descriptions of Lefschetz numbers over singular varieties. Note that for normal complex algebraic surfaces a complete answer to this problem was already given by S. Saito \cite{Saito}. It would be also an interesting problem to compare these results with the recent development in complex dynamical systems such as Abate-Bracci-Tovena \cite{A-B-T}.

In Section \ref{sec:8}, we study functorial properties of Lefschetz cycles and prove the direct and inverse image theorems for them which extend those for characteristic cycles obtained by Kashiwara-Schapira \cite{K-S}. Since the sign of the determinant $\sgn(\id-\phi^{\prime})= \pm 1$ of the linear map
\begin{equation}
\id- \phi^{\prime} \colon T_MX \longrightarrow T_MX
\end{equation}
naturally appears in the inverse image theorem, its proof is much more involved than that of the direct image theorem. To determine these very subtle signs of the determinant, the theory of currents with hyperfunction coefficients will be used. Finally, let us mention that our inverse image theorem has also an application to the explicit description of local contributions. Indeed, in Corollary \ref{cor:7-11} we prove the localizability of the global trace of $(F, \Phi)$ to the fixed point manifold $M$ without assuming any technical condition such as
\begin{equation}
\Ev_x \cap \Real_{> 1}=\emptyset \hspace{5mm} \text{for any $x\in M$}
\end{equation}
on the map $\phi$. It seems that if there exists a point $x\in M$ such that $\Ev_x \cap \Real_{> 1}\neq \emptyset$ the usual methods (see e.g. \cite[Section 9.6]{K-S} and \cite{G-M-1} etc.) for localizations do not work. Namely, our theory of Lefschetz cycles enables us to obtain the localization even when the map $\phi$ is expanding in some directions normal to $M$.

\section{Preliminary notions and results}\label{sec:2}

In this paper, we essentially follow the terminology in \cite{K-S}. For example, for a topological space $X$, we denote by $\Db(X)$ the derived category of bounded complexes of sheaves of $\Comp_X$-modules on $X$. From now on, we shall review basic notions and known results concerning Lefschetz fixed point formulas. Since we focus our attention on Lefschetz fixed point formulas for constructible sheaves in this paper, we treat here only real analytic manifolds and morphisms. Now let $X$ be a real analytic manifold. We denote by $\Dc(X)$ the full subcategory of $\Db(X)$ consisting of bounded complexes of sheaves whose cohomology sheaves are $\Real$-constructible (see \cite[Chapter VIII]{K-S} for the precise definition). Let $\phi \colon X \longrightarrow X$ be an endomorphism of the real analytic manifold $X$. Then our initial datum is a pair $(F,\Phi)$ of $F \in \Dc(X)$ and a morphism $\Phi \colon \phi^{-1}F \longrightarrow F$ in $\Dc(X)$. If the support $\supp (F)$ of $F$ is compact, $H^j(X;F)$ is a finite-dimensional vector space over $\Comp$ for any $j \in \Z$ and we can define the following important number from $(F,\Phi)$.

\begin{dfn}\label{dfn:2-1}
We set
\begin{equation}
\tr (F,\Phi) :=\dsum_{j \in \Z} (-1)^j \tr \{ H^j(X;F) \overset{\Phi}{\longrightarrow} H^j(X;F) \} \in \Comp,
\end{equation}
where the morphisms $H^j(X;F) \overset{\Phi}{\longrightarrow} H^j(X;F)$ are induced by
\begin{equation}
F \longrightarrow R\phi_*\phi^{-1}F \overset{\Phi}{\longrightarrow}R\phi_*F.
\end{equation}
We call $\tr (F,\Phi)$ the global trace of the pair $(F,\Phi)$.
\end{dfn}

Now let us set
\begin{equation}
M:= \{ x \in X \ |\ \phi(x)=x\} \subset X.
\end{equation}
This is the fixed point set of $\phi \colon X \longrightarrow X$ in $X$. Since we mainly consider the case where the fixed point set is a smooth submanifold of $X$, we use the symbol $M$ to express it. If a compact group $G$ is acting on $X$ and $\phi$ is the left action of an element of $G$, then the fixed point set is smooth by Palais's theorem \cite{Palais} (see \cite{Guillemin} for an excellent survey of this subject). Now let us consider the diagonal embedding $\delta_X \colon X \longhookrightarrow X \times X$ of $X$ and the closed embedding $h:=(\phi, \id_X) \colon X \longhookrightarrow X \times X$ associated with $\phi$. Denote by $\Delta_X$ (resp. $\Gamma_{\phi}$) the image of $X$ by $\delta_X$ (resp. $h$). Then $M \simeq \Delta_X \cap \Gamma_{\phi}$ and we obtain a chain of morphisms
\begin{eqnarray}
R\Hom_{\Comp_X}(F,F)
&\simeq & \delta_X^!(F \boxtimes \D F) \\
&\longrightarrow & \RG_{\supp(F) \cap \Delta_X} (h_*h^{-1}(F \boxtimes \D F))|_{\Delta_X} \\
&\simeq & \RG_{\supp(F) \cap \Delta_X} (h_*(\phi^{-1}F \otimes \D F))|_{\Delta_X} \\
&\overset{\Phi}{\longrightarrow}& \RG_{\supp(F) \cap \Delta_X}(h_*(F \otimes \D F))|_{\Delta_X} \\
& \longrightarrow & \RG_{\supp(F) \cap \Delta_X} (h_*\omega_X)|_{\Delta_X} \\
&\simeq & \RG_{\supp(F) \cap M} (\omega_X),
\end{eqnarray}
where $\omega_X \simeq or_X [\d X] \in \Dc(X)$ is the dualizing complex of $X$ and $\D F =R\Hom_{\Comp_X}(F, \omega_X)$ is the Verdier dual of $F$. Hence we get a morphism
\begin{equation}\label{eq:2-10}
{\rm Hom}_{\Db(X)}(F,F) \longrightarrow H^0_{\supp(F) \cap M}(X;\omega_X).
\end{equation}

\begin{dfn}[\cite{Kashiwara-2}]\label{dfn:2-2} We denote by $C(F,\Phi)$ the image of $\id_F$ by the morphism \eqref{eq:2-10} in $H_{\supp(F) \cap M}^0(X;\omega_X)$ and call it the characteristic class of $(F,\Phi)$.
\end{dfn}

\begin{thm}[\cite{Kashiwara-2}]\label{thm:2-3}
If $\supp (F)$ is compact, then the equality
\begin{equation}
\tr (F,\Phi)=\dint_X C(F,\Phi)
\end{equation}
holds. Here
\begin{equation}
\dint_X \colon H_c^n(X;or_X) \longrightarrow \Comp
\end{equation}
is the morphism induced by the integral of differential $(\d X)$-forms with compact support.
\end{thm}

Let $M= \bigsqcup_{i \in I}M_i$ be the decomposition of $M$ into connected components and
\begin{gather}
H_{\supp(F) \cap M}^0(X;\omega_X) =\bigoplus_{i \in I} H_{\supp(F) \cap M_i}^0(X;\omega_X), \\
C(F,\Phi) =\bigoplus_{i \in I} C(F,\Phi)_{M_i}
\end{gather}
the associated direct sum decomposition.

\begin{dfn}\label{dfn:2-4}
When $\supp(F) \cap M_i$ is compact, we define a complex number $c(F,\Phi)_{M_i}$ by
\begin{equation}
c(F,\Phi)_{M_i}:=\dint_X C(F,\Phi)_{M_i}
\end{equation}
and call it the local contribution of $(F,\Phi)$ from $M_i$.
\end{dfn}

By Theorem \ref{thm:2-3}, if $\supp (F)$ is compact, the global trace of $(F,\Phi)$ is the sum of local contributions:
\begin{equation}
\tr (F,\Phi)=\dsum_{i \in I} c(F,\Phi)_{M_i}.
\end{equation}
Hence one of the most important problems in the theory of Lefschetz fixed point formulas is to explicitly describe these local contributions. However the direct computation of local contributions is a very difficult task in general. Instead of local contributions, we usually consider first the following number $\tr (F|_{M_i},\Phi|_{M_i})$ which is much more easily computed. Let $M_i$ be a fixed point component such that $\supp (F) \cap M_i$ is compact.

\begin{dfn}\label{dfn:2-5}
We set
\begin{equation}
\tr (F|_{M_i},\Phi|_{M_i}):= \dsum_{j \in \Z}(-1)^j \tr\{ H^j(M_i;F|_{M_i}) \overset{\Phi|_{M_i}}{\longrightarrow} H^j(M_i;F|_{M_i})\},
\end{equation}
where the morphisms $H^j(M_i;F|_{M_i})\overset{\Phi|_{M_i}}{\longrightarrow} H^j(M_i;F|_{M_i})$ are induced by the restriction
\begin{equation}
\Phi|_{M_i} \colon F|_{M_i} \simeq (\phi^{-1}F)|_{M_i} \longrightarrow F|_{M_i}
\end{equation}
of $\Phi$.
\end{dfn}

We can easily compute this new invariant $\tr(F|_{M_i},\Phi|_{M_i})\in \Comp$ as follows. Let $M_i= \bigsqcup_{\alpha \in A}M_{i,\alpha}$ be a stratification of $M_i$ by connected subanalytic manifolds $M_{i,\alpha}$ such that $H^j(F)|_{M_{i,\alpha}}$ is a locally constant sheaf for any $\alpha \in A$ and $j \in \Z$. Namely, we assume that the stratification $M_i = \bigsqcup_{\alpha \in A} M_{i,\alpha}$ is adapted to $F|_{M_i}$.

\begin{dfn}\label{dfn:2-6}
For each $\alpha \in A$, we set
\begin{equation}
c_{\alpha}:= \dsum_{j\in \Z}(-1)^j \tr\{ H^j(F)_{x_{\alpha}} \overset{\Phi|_{\{x_{\alpha}\}}}{\longrightarrow} H^j(F)_{x_{\alpha}}\} \hspace{5mm}\in \Comp,
\end{equation}
where $x_{\alpha}$ is a reference point of $M_{i,\alpha}$.
\end{dfn}

Then we have the following very useful result due to Goresky-MacPherson.

\begin{prp}[\rm \cite{G-M-1}]\label{prp:2-7} 
In the situation as above, we have
\begin{equation}
\tr(F|_{M_i},\Phi|_{M_i})=\dsum_{\alpha \in A}c_{\alpha} \cdot \chi_c(M_{i,\alpha}),
\end{equation}
where $\chi_c$ is the Euler-Poincar{\'e} index with compact supports.
\end{prp}

In terms of the theory of topological integrals of constructible functions developed by Kashiwara-Schapira \cite{K-S} and Viro \cite{Viro} etc., we can restate this result in the following way. Since we need $\Comp$-valued constructible functions, we slightly generalize the usual notion of $\Z$-valued constructible functions.

\begin{dfn}\label{dfn:2-8}
Let $Z$ be a subanalytic set. Then we say that a $\Comp$-valued function $\varphi \colon Z \longrightarrow \Comp$ is constructible if there exists a stratification $Z= \bigsqcup_{\alpha \in A} Z_{\alpha}$ of $Z$ by subanalytic manifolds $Z_{\alpha}$ such that $\varphi|_{Z_{\alpha}}$ is a constant function for any $\alpha \in A$. We denote by $\CF(Z)_{\Comp}$ the abelian group of $\Comp$-valued constructible functions on $Z$.
\end{dfn}

Let $\varphi =\sum_{\alpha \in A} c_{\alpha} \cdot \1_{Z_{\alpha}} \in \CF(Z)_{\Comp}$ be a $\Comp$-valued constructible function with compact support on a subanalytic set $Z$, where $Z= \bigsqcup_{\alpha \in A} Z_{\alpha}$ is a stratification of $Z$ and $c_{\alpha} \in \Comp$. Then we can easily prove that the complex number $\sum_{\alpha \in A} c_{\alpha} \cdot \chi_c(Z_{\alpha})$ does not depend on the expression $\varphi=\sum_{\alpha \in A} c_{\alpha} \cdot \1_{Z_{\alpha}}$ of $\varphi$.

\begin{dfn}\label{dfn:2-9}
For a $\Comp$-valued constructible function $\varphi=\sum_{\alpha \in A}c_{ \alpha} \cdot \1_{Z_{\alpha}} \in \CF(Z)_{\Comp}$ with compact support as above, we set
\begin{equation}
\dint_Z\varphi :=\dsum_{\alpha \in A}c_{\alpha} \cdot \chi_c(Z_{\alpha}) \in \Comp
\end{equation}
and call it the topological integral of $\varphi$.
\end{dfn}

By this definition, the result of Proposition \ref{prp:2-7} can be rewritten as
\begin{equation}
\tr(F|_{M_i},\Phi|_{M_i}) =\dint_{M_i}\varphi(F,\Phi)_{M_i},
\end{equation}
where the $\Comp$-valued constructible function $\varphi(F,\Phi)_{M_i} \in \CF(M_i)_{\Comp}$ on $M_i$ is defined by
\begin{equation}\label{eq:2-23}
\varphi(F,\Phi)_{M_i}(x):=\dsum_{j \in \Z} (-1)^j\tr\{ H^j(F)_x \overset{\Phi|_{\{x\}}}{\longrightarrow} H^j(F)_x\}
\end{equation}
for $x\in M_i$.

To end this section, let us explain how the $\Comp$-valued constructible functions discussed above are related to the theory of Lagrangian cycles in \cite[Chapter IX]{K-S}. Now let $Z$ be a real analytic manifold and denote by $T^*Z$ its cotangent bundle. Recall that Kashiwara-Schapira constructed the sheaf $\L_Z$ of closed conic subanalytic Lagrangian cycles on $T^*Z$ in \cite{K-S} (in this paper, we consider Lagrangian cycles with coefficients in $\Comp$).

\begin{prp}[\cite{K-S}]\label{prp:2-10}
There exists a group isomorphism
\begin{equation}
CC \colon \CF(Z)_{\Comp} \simto \varGamma(T^*Z;\L_Z) 
\end{equation}
by which the characteristic function $\1_K$ of a closed submanifold $K \subset Z$ of $Z$ is sent to the conormal cycle $[T^*_KZ]$ in $T^*Z$.
\end{prp}
We call $CC$ the characteristic cycle map in this paper.

\section{Localization theorems and their applications}\label{sec:3}

Let $X$, $\phi \colon X \longrightarrow X$, $M= \bigsqcup_{i \in I}M_i$, $F \in \Dc(X)$, $\Phi \colon \phi^{-1}F \longrightarrow F$ etc. be as in Section \ref{sec:2}. In this section, we fix a fixed point component $M_i$ and always assume that $\supp(F) \cap M_i$ is compact.

\begin{dfn}\label{dfn:3-1}
We say that the global trace $\tr(F,\Phi)$ is localizable to $M_i$ if the equality
\begin{equation}
c(F,\Phi)_{M_i}=\tr(F|_{M_i},\Phi|_{M_i})
\end{equation}
holds.
\end{dfn}

By Proposition \ref{prp:2-7}, once the global trace is localizable to $M_i$, the local contribution $c(F,\Phi)_{M_i}$ of $(F, \Phi)$ from $M_i$ can be very easily computed. Since we always consider the same fixed point component $M_i$ in this section, we denote $M_i$, $c(F,\Phi)_{M_i}$ etc. simply by $M$, $c(F,\Phi)_M$ etc. respectively. From now on, we shall give a useful criterion for the localizability of the global trace to $M$. First let us consider the natural morphism
\begin{equation}
\phi^{\prime} \colon T_{M_{\reg}}X \longrightarrow T_{M_{\reg}}X
\end{equation}
induced by $\phi \colon X \longrightarrow X$, where $M_{\reg}$ denotes the set of regular points in $M$. Since $M_{\reg}$ is not always connected in the real analytic case, the rank of $T_{M_{\reg}}X$ may vary depending on the connected components of $M_{\reg}$.

\begin{dfn}\label{dfn:3-2}
For $x \in M_{\reg}$, we set
\begin{equation}
\Ev_x :=\{ \text{the eignevalues of $\phi^{\prime}_x \colon (T_{M_{\reg}}X)_x \longrightarrow (T_{M_{\reg}}X)_x$} \} \subset \Comp.
\end{equation}
\end{dfn}

We also need the specialization functor
\begin{equation}
\nu_{M_{\reg}} \colon \Db(X) \longrightarrow \Db(T_{M_{\reg}}X)
\end{equation}
along $M_{\reg} \subset X$. In order to recall the construction of this functor, consider the standard commutative diagram:
\begin{equation}
\xymatrix{
T_{M_{\reg}}X \ar[d]^{\tau} \ar@{^{(}->}[r]^s & \tl{X_{M_{\reg}}} \ar[d]^p & \Omega_X \ar@{_{(}->}[l]_j \ar[dl]^{\tl{p}} \\
M_{\reg} \ar@{^{(}->}[r]^i& X,}
\end{equation}
where $\tl{X_{M_{\reg}}}$ is the normal deformation of $X$ along $M_{\reg}$ and $t \colon \tl{X_{M_{\reg}}} \longrightarrow \Real$ is the deformation parameter. Recall that $\Omega_X$ is defined by $t > 0$ in $\tl{X_{M_{\reg}}}$. Then the specialization $\nu_{M_{\reg}}(F)$ of $F$ along $M_{\reg}$ is defined by
\begin{equation}
\nu_{M_{\reg}}(F):=s^{-1}Rj_*\tl{p}^{-1}(F).
\end{equation}
Note that $\nu_{M_{\reg}}(F)$ is a conic object in $\Db(T_{M_{\reg}}X)$ whose support is contained in the normal cone $C_{M_{\reg}}(\supp(F))$ to $\supp(F)$ along $M_{\reg}$. Since $F$ is $\Real$-constructible, $\nu_{M_{\reg}}(F)$ is also $\Real$-constructible. By construction, there exists a natural morphism
\begin{equation}
\Phi^{\prime} \colon (\phi^{\prime})^{-1}\nu_{M_{\reg}}(F) \longrightarrow \nu_{M_{\reg}}(F)
\end{equation}
induced by $\Phi \colon \phi^{-1}F \longrightarrow F$. In the sequel, let us assume the conditions:
\begin{enumerate}
\item $\supp(F) \cap M$ is compact and contained in $M_{\reg}$.
\item $1 \notin \Ev_x$ for any $x \in \supp(F) \cap M_{\reg}$.
\end{enumerate}
The condition (ii) implies that the graph of $\phi$ in $X\times X$ intersects cleanly (see \cite[Definition 4.1.5]{K-S}) with the diagonal set $\Delta_X \simeq X$ in an open neighborhood of $\supp(F) \cap M_{\reg}$. It follows also from the condition (ii) that for an open neighborhood $U$ of $\supp(F) \cap M_{\reg}$ in $M_{\reg}$ the fixed point set of $\phi^{\prime}|_{\tau^{-1}(U)} \colon \tau^{-1}(U) \longrightarrow \tau^{-1}(U)$ is contained in the zero-section $M_{\reg}$ of $T_{M_{\reg}}X$. Set $\tl{U}=\tau^{-1}(U)$, $\tl{F}=\nu_{M_{\reg}}(F)|_{\tl{U}}$ and $\tl{\Phi} =\Phi^{\prime}|_{\tl{U}} \colon (\phi^{\prime}|_{\tl{U}})^{-1}\tl{F} \longrightarrow \tl{F}$. Then also for the pair $(\tl{F},\tl{\Phi})$, we can define the characteristic class $C(\tl{F},\tl{\Phi})\in H_{\supp(F) \cap M_{\reg}}^0(\tl{U};\omega_{\tl{U}})$.

\begin{prp}\label{prp:3-3}
In the situation as above, the local contribution $c(F,\Phi)_M$ from $M$ is equal to $\dint_{\tl{U}}C(\tl{F},\tl{\Phi})$.
\end{prp}

\begin{proof}
The proof is similar to that of \cite[Proposition 9.6.11]{K-S}. Since the construction of the characteristic class $C(F,\Phi)_M \in H_{\supp(F) \cap M}^0(X;\omega_X)$ is local around $\supp(F) \cap M$ (see \cite[Remark 9.6.7]{K-S}) and $X \setminus (M \setminus U)$ is invariant by $\phi$, we may replace $X$, $M$ etc. by $X\setminus (M \setminus U)$, $U$ etc. respectively. Then the proof follows from the commutativity of the diagram \eqref{diag:3-6} below. Here we denote $T_MX$ simply by $\Vb$ and the morphism $\tl{h} \colon T_MX \longrightarrow T_MX \times T_MX$ is defined by $\tl{h}=(\phi^{\prime},\id)$. We also used the natural isomorphism $\D \nu_M(F) \simeq \nu_M(\D F)$. Let us explain the construction of the morphism ${\bf A}$ in the diagram \eqref{diag:3-6}. Consider the commutative diagram:

\begin{equation}
\xymatrix@R=3mm@C=5mm{
T_{M\times M}(X\times X) \ar@{^{(}->}[rr]^{s_1} & & \tl{(X\times X)_{M\times M}} & & \Omega_{X\times X} \ar@{_{(}->}[ll]_{j_1} \ar@{->>}[rr]^{\tl{p_1}}& & X\times X \\
 & \Box & & \Box & & \Box & \\
T_MX \ar@{^{(}->}[rr]^{s} \ar@{^{(}->}[uu]^{\delta_{T_MX}} & & \tl{X_M} \ar@{^{(}->}[uu]^{\tl{\delta^{\prime}}} & & \Omega_X \ar@{_{(}->}[ll]_{j} \ar@{->>}[rr]^{\tl{p}} \ar@{^{(}->}[uu]^{\tl{\delta}} & & X, \ar@{^{(}->}[uu]^{\delta_X}}
\end{equation}
where $\tl{(X\times X)_{M \times M}}$ is the normal deformation of $X\times X$ along $M\times M$ and \\$t_1 \colon \tl{(X\times X)_{M\times M}} \longrightarrow \Real$ is the deformation parameter such that $\Omega_{X\times X}$ is defined by $t_1 > 0$ in $\tl{(X\times X)_{M \times M}}$. Then the morphism ${\bf A}$ is constructed by the morphisms of functors
\begin{eqnarray}
\delta_X^!
&\longrightarrow & \delta_X^!R\tl{p_1}_*\tl{p_1}^{-1} \\
&\simeq & R\tl{p}_* \tl{\delta}^!\tl{p_1}^{-1} \\
&\simeq & Rp_* \tl{\delta^{\prime}}^! Rj_{1*}\tl{p_1}^{-1} \\
&\longrightarrow & Rp_* \tl{\delta^{\prime}}^! s_{1*}s_1^{-1}Rj_{1*}\tl{p_1}^{-1} \\
&\simeq & Rp_*s_* \delta_{T_MX}^!s_1^{-1}Rj_{1*}\tl{p_1}^{-1}.
\end{eqnarray}
The other horizontal arrows in the diagram \eqref{diag:3-6} are constructed similarly.
{\footnotesize
\begin{equation}\label{diag:3-6}
\xymatrix@R=8mm@C=4mm{
R{\rm Hom}(F,F) \ar[rr] & & R{\rm Hom}(\nu_M(F),\nu_M(F)) \\
\RG_{\Delta_X}(X\times X;F \boxtimes \D F) \ar[u]^{\wr} \ar[r]^{\hspace*{-5mm}{\bf A}} \ar[dd] & \RG_{\Delta_{\Vb}}(\Vb\times \Vb;\nu_{M\times M}(F \boxtimes \D F)) \ar[d]& \RG_{\Delta_{\Vb}}(\Vb\times \Vb;\nu_M(F) \boxtimes \D \nu_M(F)) \ar[l] \ar[d] \ar[u]^{\wr}\\
 & \RG_{M}(\Vb;\tl{h}^{-1}\nu_{M\times M}(F \boxtimes \D F)) \ar[d] & \RG_{M}(\Vb;\phi^{\prime -1}\nu_M(F) \otimes \D \nu_M(F)) \ar[l] \ar[d] \\
\RG_M(X; \phi^{-1}F \otimes \D F) \ar[d]^{\Phi} \ar[r]^{\sim} & \RG_M(\Vb; \nu_M(\phi^{-1}F \otimes \D F)) \ar[d]^{\Phi} & \RG_M(\Vb;\nu_M(\phi^{-1}F) \otimes \D \nu_M(F)) \ar[l] \ar[d]^{\Phi} \\
\RG_M(X;F\otimes \D F) \ar[d] \ar[r]^{\sim} & \RG_M(\Vb;\nu_M(F\otimes \D F)) \ar[d] & \RG_M(\Vb;\nu_M(F) \otimes \D \nu_M(F)) \ar[l] \ar[d] \\
\RG_M(X;\omega_X) \ar[d] \ar[r]^{\sim} & \RG_M(\Vb;\nu_M(\omega_X)) & \RG_M(\Vb;\omega_{\Vb}) \ar[d] \ar@{-}[l]_{\sim}\\
\Comp \ar@{=}[rr]& & \Comp. }
\end{equation}}\qed
\end{proof}

\begin{thm}\label{thm:3-4}
In the situation as above, assume moreover that
\begin{equation}\label{eq:3-6}
\Ev_x \cap \Real_{\geq 1} =\emptyset
\end{equation}
for any $x \in \supp(F) \cap M \subset M_{\reg}$. Then the localization
\begin{equation}
c(F,\Phi)_M =\tr(F|_M,\Phi|_M)
\end{equation}
holds.
\end{thm}

\begin{proof}
The proof is similar to that of \cite[Proposition 9.6.12]{K-S}. Since
\begin{equation}
\supp(F) \cap M \subset M_{\reg}
\end{equation}
is compact, there exists an open neighborhood $U$ of $\supp(F) \cap M$ in $M_{\reg}$ such that \eqref{eq:3-6} holds for any $x \in U$.

As in the proof of Proposition \ref{prp:3-3}, we may replace $X$, $M$ etc. by $X\setminus (M \setminus U)$, $U$ etc. respectively. By the homotopy invariance of local contributions (\cite[Proposition 9.6.8]{K-S}), replacing $\phi^{\prime}$ by $\lambda \phi^{\prime}$ for $0<\lambda<1$ does not affect
\begin{gather}
C(\nu_M(F),\Phi^{\prime}) \in H_{\supp(F) \cap M}^0(T_MX;\omega_{T_MX}), \\
\tr(\nu_M(F)|_M,\Phi^{\prime}|_M)=\tr(F|_M,\Phi|_M).
\end{gather}
Since $\supp(F) \cap M$ is compact, we may take sufficiently small $0 <\lambda \ll1$ so that the condition
\begin{equation}
\Ev_x \subset \{ z \in \Comp\ |\ |z|<1\} \hspace{5mm}\text{for any $x \in U$}
\end{equation}
is satisfied. Then just in the same way as in Step (a) of the proof of \cite[Proposition 9.6.12]{K-S}, we can prove
\begin{eqnarray}
\dint_{T_MX}C(\nu_M(F),\Phi^{\prime})
&=& \tr(\nu_M(F)|_M,\Phi^{\prime}|_M) \\
&=& \tr(F|_M,\Phi|_M).
\end{eqnarray}
Since we have
\begin{equation}
c(F,\Phi)_M=\dint_{T_MX}C(\nu_M(F),\Phi^{\prime})
\end{equation}
by Proposition \ref{prp:3-3}, this completes the proof.
\qed
\end{proof}

Similarly, in the complex case we have the following.

\begin{thm}\label{thm:3-5}
In the situation as above, assume moreover that $X$ and $\phi \colon X \longrightarrow X$ are complex analytic and $F \in \Db_c(X)$ i.e. $F$ is $\Comp$-constructible. Assume also that $M$ or $\supp(F) \cap M \subset M_{\reg}$ is smooth and compact (this condition will be removed in Section \ref{sec:6}). Then the localization
\begin{equation}
c(F,\Phi)_M=\tr(F|_M,\Phi|_M)
\end{equation}
holds.
\end{thm}

\begin{proof}
By our assumptions, the fixed point component $M$ of $\phi$ is a complex analytic subset of $X$ and $T_{M_{\reg}}X$ is a holomorphic vector bundle over $M_{\reg}$. Let $\supp(F) \cap M= \bigsqcup_{j \in J}V_j$ be the decomposition of $\supp(F) \cap M$ into connected components. Since the set $\Ev_x \subset \Comp$ of eigenvalues depends holomorphically on $x \in M_{\reg}$, $\Ev_x$ is constant on each connected component $V_j$. Hence, by the $\Comp^{\times}$-conicness of $\nu_{M_{\reg}}(F)\in \Db_c(T_{M_{\reg}}X)$ and the homotopy invariance of local contributions, for each $j \in J$ we may replace $\phi^{\prime}$ by $\lambda \phi^{\prime}$ for $|\lambda-1|\ll1$ on an open neighborhood of $\tau^{-1}(V_j) \subset T_{M_{\reg}}X$ and assume that
\begin{equation}
\Ev_x \subset \{ z \in \Comp \ |\ |z| \neq 1, \ \ z \notin \Real\} \sqcup \{0\}
\end{equation}
for any $x \in V_j$.

Then with the help of Proposition \ref{prp:3-3} and the arguments in the proof of \cite[Proposition 9.6.12]{K-S}, we may argue as in the proof of \cite[Proposition 9.6.12]{K-S}and \cite[Corollary 9.6.16]{K-S}.
\qed
\end{proof}

\begin{rem}\label{rem:3-6}
In Section \ref{sec:6}, we will generalize Theorem \ref{thm:3-5} to the case where $M$ nor $\supp(F) \cap M$ is smooth. To treat this more general case where the set $\Ev_x$ may vary depending on $x \in \supp(F) \cap M$, we need some precise arguments on Lefschetz cycles which will be introduced in the next section. One naive idea to treat this case would be to cover $\supp(F) \cap M$ by sufficiently small closed subsets $Z_i \subset \supp(F) \cap M$ and use the local contributions of $( \nu_{M_{\reg}}(F) )_{\tau^{-1}Z_i}$ to compute that of $\nu_{M_{\reg}}(F)$ by a Mayer-Vietoris type argument. However this very simple idea does not work, because we cannot apply \cite[Proposition 9.6.2]{K-S} to constructible sheaves with ``non-compact" support such as $( \nu_{M_{\reg}}(F) )_{\tau^{-1}Z_i}$ to justify the Mayer-Vietoris type argument.
\end{rem}

\begin{cor}\label{cor:3-7}
Let $X$ be a complex manifold, $\phi \colon X \longrightarrow X$ a holomorphic map and $V \subset X$ a $\phi$-invariant compact analytic subset. Assume that the fixed point set $M=\{x \in X \ |\  \phi(x)=x\} \subset X$ of $\phi$ satisfies the following conditions. 
\begin{enumerate}
\item $V \cap M \subset M_{\reg}$,
\item $1 \notin \Ev_x$ for any $x \in V \cap M$,
\item $M$ or $V \cap M$ is smooth and compact (this condition will be removed in Section \ref{sec:6}).
\end{enumerate}
Then we have
\begin{equation}
\dsum_{j \in \Z} (-1)^j \tr\{ H^j(V;\Comp_V) \overset{(\phi|_V)^*}{\longrightarrow}H^j(V;\Comp_V) \}=\chi(V \cap M).
\end{equation}
\end{cor}

\begin{proof}
Set $F=\Comp_V\in \Db_c(X)$ and let $\Phi \colon \phi^{-1}F \longrightarrow F$ be the natural morphism $\phi^{-1}\Comp_V\simeq \Comp_{\phi^{-1}(V)} \longrightarrow \Comp_V$. Then we have
\begin{equation}
\dsum_{j\in \Z}(-1)^j\tr\{H^j(X;F) \overset{\Phi}{\longrightarrow}H^j(X;F)\}= \dsum_{j \in \Z}(-1)^j\tr\{H^j(V;\Comp_V) \overset{(\phi|_V)^*}{\longrightarrow} H^j(V;\Comp_V) \}.
\end{equation}
By Theorem \ref{thm:3-5}, this number is equal to $\tr(F|_M,\Phi|_M)=\dint_M \varphi(F,\Phi)_M$, where \\$\varphi(F,\Phi)_M \colon M \longrightarrow \Comp$ is the $\Comp$-valued constructible function on $M$ defined as in \eqref{eq:2-23}. Since $\phi|_M$ is the identity map of $M$, $\varphi(F,\Phi)_M=\1_{V \cap M}$ and the result follows.
\qed
\end{proof}

\begin{exm}\label{exm:3-8}
Let $G_n=SL_n(\Comp)$ and let $B_n \subset G_n$ be the Borel subgroup of $G_n$ consisting of upper triangular matrices. Then the homogeneous space $X=G_n/B_n$ is a flag manifold. Take an element
\begin{equation}
g={\rm diag}(\underbrace{\lambda_1,\cdots,\lambda_1}_{\text{$n_1$-times}}, \underbrace{\lambda_2,\cdots,\lambda_2}_{\text{$n_2$-times}}, \cdots, \underbrace{\lambda_k,\cdots,\lambda_k}_{\text{$n_k$-times}})
\end{equation}
$(n=n_1+\cdots +n_k)$ in $B_n$ such that $\lambda_i \neq \lambda_j$ for any $i \neq j$, where ${\rm diag}(\cdots)$ denotes a diagonal matrix. Let $\phi \colon X \longrightarrow X$ be the left action $l_g \colon X \simto X$ by $g \in B_n \subset G_n$. Then it is easy to see that the fixed point set $M$ of $\phi$ is a smooth complex submanifold of $X$. More precisely, $M$ is isomorphic to the disjoint union of $\dfrac{n!}{n_1!\cdots n_k!}$ copies of the product of flag manifolds
\begin{equation}
G_{n_1}/B_{n_1} \times \cdots \times G_{n_k}/B_{n_k}.
\end{equation}
Therefore the assumptions of Corollary \ref{cor:3-7} are satisfied for any $\phi$-invariant analytic subset $V$ of $X$, if $1 \notin \Ev_x$ for any $x \in M$ (we expect this is always true). Since $g \in B_n$, as a $\phi$-invariant analytic subset $V$ we can take any Schubert variety in $X$.
\end{exm}

\begin{exm}\label{exm:3-9}
Let us consider a special case of Example \ref{exm:3-8} above. Let $X=G_3/B_3$ be the flag manifold consisting of full flags in $\Comp^3$ and $\phi =l_g \colon X \simto X$ the left action by the element
\begin{equation}
g=\begin{pmatrix} \alpha & 0 & 0 \\ 0 & \alpha & 0 \\ 0 & 0 & \beta \end{pmatrix} \in B_3 \subset G_3,
\end{equation}
where $\alpha \neq \beta$ are non-zero complex numbers. In this case, the fixed point set $M \subset X$ of $\phi$ is the disjoint union of $3$ copies of $\Comp\P^1$'s. Let $X= \bigsqcup_{\sigma \in \S_3} B_3\sigma B_3= \bigsqcup_{\sigma \in \S_3}X_{\sigma}$ be the Bruhat decomposition of $X=G_3/B_3$. Here an element $\sigma$ of the symmetric group $\S_3$ is identified with the matrix $(\delta_{i, \sigma(j)})_{1\leq i,j \leq 3}\in G_3$ (see e.g. \cite{H-T-T} for the detail of this subject), where $\delta_{ij}$ is Kronecker's delta. In this case, $X_{(1,3)}$ is the unique open dense Schubert cell in $X$. Set $V=X\setminus X_{(1,3)} = \bigsqcup_{\sigma \neq (1,3)} X_{\sigma}$. Then $V$ is a $\phi$-invariant analytic subset of $X$ and we can check that the assumptions of Corollary \ref{cor:3-7} are satisfied. For example, let $N_-$ be a unipotent subgroup of $B_3$ defined by
\begin{equation}
N_- =\left\{ \left. \begin{pmatrix} 1 & 0& 0 \\ a & 1 & 0 \\ b & c & 1\end{pmatrix} \ \right| \ a,b,c \in \Comp \right\}
\end{equation}
and consider the open embedding
\begin{equation}
\begin{array}{ccccc}
\Comp^3 & \simeq & N_- & \longhookrightarrow & X=G_3/B_3. \\
\inun & & \inun & & \inun \\
\begin{pmatrix} a \\ b \\ c \end{pmatrix} & \longmapsto & \begin{pmatrix} 1 & 0& 0 \\ a & 1 & 0 \\ b & c & 1 \end{pmatrix} & \longmapsto & \begin{pmatrix} 1 & 0& 0 \\ a & 1 & 0 \\ b & c & 1 \end{pmatrix}B_3.
\end{array}
\end{equation}
We denote the image of this open embedding by $U$ and identify it with $\Comp^3=\{(a,b,c) \in \Comp^3\}$. Then in $U\simeq \Comp^3$ the fixed point set $M=\Comp\P^1 \sqcup \Comp\P^1 \sqcup \Comp\P^1$ of $\phi$ is given by
\begin{equation}
M \cap U=\{ (a,b,c) \in \Comp^3 \ |\ b=c=0\}
\end{equation}
and $1 \notin \Ev_x =\left\{ \dfrac{\beta}{\alpha}, \dfrac{\beta}{\alpha}\right\}$ for any $x \in M \cap U$. The assumptions of Corollary \ref{cor:3-7} are thus verified. Moreover in $U \simeq \Comp^3$ the $\phi$-invariant analytic subset $V =X\setminus X_{(1,3)}$ is given by
\begin{equation}
V \cap U =\{ (a,b,c)\in \Comp^3 \ |\ b(b-ac)=0\}.
\end{equation}
From this, we see that the fixed point component $\overline{M \cap U} \simeq \Comp\P^1$ of $M$ is totally embedded in the singular set of $V$.
\end{exm}

\section{Definition of Lefschetz cycles}\label{sec:4}

In this section, we construct certain Lagrangian cycles which encode the local contributions discussed in previous sections into topological objects. We inherit the notations in Section \ref{sec:2} and \ref{sec:3}. Now assume that the fixed point set $M=\{ x \in X \ |\ \phi(x)=x\}$ of $\phi \colon X \longrightarrow X$ is a submanifold of $X$. However here we do not assume that $M$ is connected. We also assume that $\Delta_X$ intersects with $\Gamma_{\phi}=\{ (\phi(x),x) \in X \times X \ |\ x \in X\}$ cleanly along $M$ in $X \times X$. Identifying $\Gamma_{\phi}$ with $X$ by the second projection $X \times X \longrightarrow X$, we obtain a natural identification $M=\Gamma_{\phi} \cap \Delta_X$. We also identify $T_{\Delta_X}(X \times X)$ (resp. $T^*_{\Delta_X}(X\times X)$) with $TX$ (resp. $T^*X$) by the first projection $T(X\times X) \simeq TX \times TX \longrightarrow TX$ (resp. $T^*(X\times X) \simeq T^*X \times T^*X \longrightarrow T^*X$) as usual. Then, under the above assumptions, we see that the natural morphism
\begin{equation}
T_M\Gamma_{\phi} \longrightarrow T_{\Delta_X}(X\times X) \simeq TX
\end{equation}
induced by the inclusion map $\Gamma_{\phi} \longhookrightarrow X \times X$ is injective. Hence the image of this morphism is a subbundle of $M \utimes{X}TX$ (whose rank may vary depending on the connected components of $M$). 

\begin{dfn}
We denote by $\tLb$ the subbundle of $M\utimes{X}TX$ constructed above.
\end{dfn}

The following lemma will be obvious.

\begin{lem}
The subset $T^*_{\Gamma_{\phi}}(X\times X) \cap T^*_{\Delta_X}(X\times X)$ of $(\Gamma_{\phi} \cap \Delta_X) \utimes{\Delta_X} T^*_{\Delta_X}(X\times X) \simeq M \utimes{X}T^*X$ consists of covectors which are orthogonal to $\tLb \subset M\utimes{X}TX$ by the natural perfect pairing $(M\utimes{X}TX) \times_M (M \utimes{X}T^*X) \longrightarrow \Real$.
\end{lem}

By this lemma, we see that $T^*_{\Gamma_{\phi}}(X\times X) \cap T^*_{\Delta_X}(X\times X)$ is a subbundle of $M\utimes{X}T^*X$ (whose rank may vary depending on the connected components of $M$).

\begin{dfn}
We denote the subbundle $T^*_{\Gamma_{\phi}}(X\times X) \cap T^*_{\Delta_X}(X\times X)$ of $M\utimes{X}T^*X$ by $\Lb$ and call it the Lefschetz bundle associated with $\phi \colon X \longrightarrow X$.
\end{dfn}

The Lefschetz bundles satisfy the following nice property.

\begin{prp}\label{prp:4-4}
The natural surjective morphism $\rho \colon M\utimes{X}T^*X \longtwoheadrightarrow T^*M$ induces an isomorphism $\Lb \simto T^*M$.
\end{prp}

\begin{proof}
Since the rank of the Lefschetz bundle $\Lb$ is locally constant and equal to that of $T^*M$, it suffices to show that the morphism
\begin{equation}
\rho|_{\Lb} \colon \Lb \longrightarrow T^*M
\end{equation}
is injective. This follows immediately from the equality
\begin{equation}
TM +\tLb =M\utimes{X}TX
\end{equation}
obtained by our hypothesis.
\qed
\end{proof}

\noindent From now on, by Proposition \ref{prp:4-4} we shall identify the Lefschetz bundle $\Lb$ with $T^*M$.

Now let $F$ be an object of $\Dc(X)$ and $\Phi \colon \phi^{-1}F \longrightarrow F$ a morphism in $\Dc(X)$. To these data $(F, \Phi)$, we can associate a conic Lagrangian cycle in the Lefschetz bundle $\Lb \simeq T^*M$ as follows. Denote by $\pi_X \colon T^*X \longrightarrow X$ the natural projection and recall that we have the functor
\begin{equation}
\mu_{\Delta_X} \colon \Db(X\times X) \longrightarrow \Db(T^*_{\Delta_X}(X\times X))
\end{equation}
of microlocalization which satisfies
\begin{equation}
R\pi_{X*}\mu_{\Delta_X} \simeq \delta_X^! \simeq \delta_X^{-1}\RG_{\Delta_X}.
\end{equation}
Recall also that the micro-support $\SS(F)$ of $F$ is a closed conic subanalytic Lagrangian subset of $T^*X$ and the support of $\mu_{\Delta_X}(F\boxtimes \D F)$ is contained in $\SS(F) \subset T^*X \simeq T^*_{\Delta_X}(X\times X)$. Then we have a chain of natural morphisms:
\begin{eqnarray}
R{\rm Hom}_{\Comp_X}(F,F)
&\simeq & \RG(X;\delta_X^!(F \boxtimes \D F)) \label{eq:4-6}\\
&\simeq & \RG_{\SS(F)}(T^*X;\mu_{\Delta_X}(F \boxtimes \D F)) \\
&\longrightarrow & \RG_{\SS(F)}(T^*X; \mu_{\Delta_X}(h_*h^{-1}(F \boxtimes \D F))) \\
&\simeq & \RG_{\SS(F)}(T^*X;\mu_{\Delta_X}(h_*(\phi^{-1}F \otimes \D F))) \\
&\overset{\Phi}{\longrightarrow} & \RG_{\SS(F)}(T^*X; \mu_{\Delta_X}(h_*(F \otimes \D F))) \\
& \longrightarrow & \RG_{\SS(F)}(T^*X;\mu_{\Delta_X}(h_*\omega_X)).\label{eq:4-11}
\end{eqnarray}

\begin{lem}\label{lem:4-5}
\begin{enumerate}
\item The support of $\mu_{\Delta_X}(h_*\omega_X)$ is contained in $\Lb$.
\item The restriction of $\mu_{\Delta_X}(h_*\omega_X)$ to $\Lb\simeq T^*M$ is isomorphic to $\pi_M^{-1}\omega_M$, where \\ $\pi_M \colon T^*M \longrightarrow M$ is the natural projection.
\end{enumerate}
\end{lem}

\begin{proof}
(i)\ \ By $\SS(h_*\omega_X) =T^*_{\Gamma_{\phi}}(X\times X)$, we obtain
\begin{equation}
\supp (\mu_{\Delta_X}(h_*\omega_X)) \subset T^*_{\Gamma_{\phi}}(X\times X) \cap T^*_{\Delta_X}(X\times X) =\Lb.
\end{equation}

(ii)\ \ Let $i_M \colon M \longhookrightarrow X$ be the inclusion map. Since we have
\begin{eqnarray}
\RG_{\Delta_X}(h_*\omega_X) &\simeq & \delta_{X*}i_{M*}\omega_M, \\
\mu_{\Delta_X}(\delta_{X*}i_{M*}\omega_M) &\simeq & \pi_X^{-1}i_{M*}\omega_M,
\end{eqnarray}
we obtain a morphism
\begin{eqnarray}
\pi_X^{-1}i_{M*}\omega_M
&\simeq & \mu_{\Delta_X}(\RG_{\Delta_X}(h_*\omega_X)) \\
&\longrightarrow & \mu_{\Delta_X}(h_*\omega_X).
\end{eqnarray}
It remains to show that the restriction of this morphism to $\Lb \subset T^*X$ is an isomorphism. Let $p \in \Lb \subset T^*X \simeq T^*_{\Delta_X}(X\times X)$ be a point. Then we have
\begin{equation}
H^j(\mu_{\Delta_X}(h_*\omega_X))_p \simeq \injlim{U,Z} H^j_Z(U;h_*\omega_X)
\end{equation}
for any $j\in \Z$, where $U$ (resp. $Z$) ranges through open (resp. closed) subsets of $X\times X$ such that the point $\pi_X(p) \in \Delta_X$ is contained in $U$ (resp. the normal cone $C_{\Delta_X}(Z)_{\pi_X(p)} \subset (TX)_{\pi_X(p)}$ is contained in $\{ v \in (TX)_{\pi_X(p)} \ |\ \langle v, p \rangle >0\} \cup \{0\}$).

Since $\Delta_X$ intersects with $\Gamma_{\phi}$ cleanly along $M=\Delta_X \cap \Gamma_{\phi}$, for the closed subsets $Z \subset X \times X$ above we have
\begin{equation}
Z \cap \supp (h_*\omega_X) =Z \cap \Gamma_{\phi}=M=\Delta_X \cap \supp(h_*\omega_X)
\end{equation}
in an open neighborhood of $\pi_X(p)\in \Delta_M \simeq M \subset \Delta_X \simeq X$. Namely, for sufficiently small $U$ we have an isomorphism
\begin{equation}
\RG_Z(U;\RG_{\Delta_X}(h_*\omega_X)) \overset{\sim}{\longrightarrow} \RG_Z(U;h_*\omega_X).
\end{equation}
This implies that the morphism
\begin{equation}
\mu_{\Delta_X}(\RG_{\Delta_X}(h_*\omega_X))_p \longrightarrow \mu_{\Delta_X}(h_*\omega_X)_p
\end{equation}
is an isomorphism.
\qed
\end{proof}

\begin{rem}
By taking the Fourier-Sato transform of Proposition \ref{prp:7-1} below, we can obtain a more functorial proof of Lemma \ref{lem:4-5}. However here we gave another proof in order to look at the structure of $\mu_{\Delta_X}(h_*\omega_X)$ more directly.
\end{rem}

By Lemma \ref{lem:4-5} there exists an isomorphism
\begin{equation}\label{eq:4-21}
\mu_{\Delta_X}(h_*\omega_X) \simeq (i_{\Lb})_*\pi_M^{-1}\omega_M,
\end{equation}
where $i_{\Lb} \colon \Lb \longhookrightarrow T^*X$ is the inclusion map. In what follows, we sometimes omit the symbol $(i_{\Lb})_*$ in the above identification \eqref{eq:4-21}. Combining the chain of morphisms \eqref{eq:4-6}-\eqref{eq:4-11} with the isomorphism \eqref{eq:4-21}, we obtain a morphism
\begin{equation}\label{eq:4-22}
{\rm Hom}_{\Db(X)}(F,F) \longrightarrow H^0_{\SS(F) \cap \Lb}( \Lb; \pi_M^{-1}\omega_M).
\end{equation}

\begin{dfn}
We denote by $LC(F,\Phi)$ the image of $\id_F \in {\rm Hom}_{\Db(X)}(F,F)$ in \\$H^0_{\SS(F) \cap \Lb}(\Lb; \pi_M^{-1}\omega_M)$ by the morphism \eqref{eq:4-22}.
\end{dfn}

\begin{lem}
$\SS(F) \cap \Lb$ is contained in a closed conic subanalytic Lagrangian subset of $\Lb \simeq T^*M$. Here we regard $\Lb$ as a symplectic manifold by using the standard symplectic structure of $T^*M$.
\end{lem}

\begin{proof}
Recall that $\Lb \subset M\utimes{X}T^*M$ is identified with $T^*M$ by the morphism $\rho \colon M \utimes{X}T^*X \longtwoheadrightarrow T^*M$. Now let $X= \bigsqcup_{\alpha \in A} X_{\alpha}$ be a $\mu$-stratification (see \cite[Chapter VIII]{K-S} for the definition) adapted to $F \in \Dc(X)$ and the closed subset $M \subset X$. Then there exists a subset $B \subset A$ such that $M= \bigsqcup_{\beta \in B}X_{\beta}$ and $\SS(F) \subset \bigsqcup_{\alpha \in A}T^*_{X_{\alpha}}X$. Since $\Lb$ is contained in $M \utimes{X}T^*X$, we have
\begin{equation}
\SS(F) \cap \Lb \subset \( \bigsqcup_{\beta \in B}T^*_{X_{\beta}}X \) \cap \Lb.
\end{equation}
But for any strata $X_{\beta}$ contained in $M$ ($\Longleftrightarrow \beta \in B$) the isomorphism
\begin{equation}
\rho|_{\Lb} \colon \Lb \simto T^*M
\end{equation}
induces an isomorphism
\begin{equation}
T^*_{X_{\beta}}X \cap \Lb \simto T^*_{X_{\beta}}M.
\end{equation}
Hence, via the identification $\Lb \simto T^*M$, we obtain the inclusion
\begin{equation}
\SS(F) \cap \Lb \subset \bigsqcup_{\beta \in B}T^*_{X_{\beta}}M.
\end{equation}
\qed
\end{proof}

\begin{dfn}
Choose a closed conic subanalytic Lagrangian subset $\Lambda$ of $\Lb \simeq T^*M$ such that $\SS(F) \cap \Lb \subset \Lambda$. We consider $LC(F,\Phi)$ as an element of $H^0_{\Lambda}(\Lb;\pi_M^{-1}\omega_M)$ and call it the Lefschetz cycle associated with the pair $(F,\Phi)$.
\end{dfn}

As a basic property of Lefschetz cycles, we have the following homotopy invariance. Let $I=[0,1]$ and let $\phi \colon X\times I \longrightarrow X$ be the restriction of a morphism of real analytic manifolds $X \times \Real \longrightarrow X$. For $t \in I$, let $i_t \colon X \longhookrightarrow X \times I$ be the injection defined by $x \longmapsto (x,t)$ and set $\phi_t :=\phi \circ i_t \colon X \longrightarrow X$. Assume that the fixed point set of $\phi_t$ in $X$ is smooth and does not depend on $t \in I$. We denote this fixed point set by $M$. Let $F\in \Dc(X)$ and consider a morphism $\Phi \colon \phi^{-1}F \longrightarrow p^{-1}F$ in $\Dc(X \times I)$, where $p \colon X \times I \longrightarrow X$ is the projection. We set
\begin{equation}
\Phi_t :=\Phi|_{X\times \{t\}} \colon \phi_t^{-1}F \longrightarrow F
\end{equation}
for $t \in I$. We denote the Lefschetz bundle associated with $\phi_t$ by $\Lb_t \simeq T^*M$.

\begin{prp}\label{prp:4-9}
Assume that $\supp(F) \cap M$ is compact and $\Lb_t$ does not depend on $t \in I$. Then the Lefschetz cycle $LC(F,\Phi_t) \in H^0_{\SS(F) \cap T^*M}(T^*M; \pi_M^{-1}\omega_M)$ does not depend on $t \in I$.
\end{prp}

\begin{proof}
The proof proceeds completely in the same way as that of \cite[Proposition 9.6.8]{K-S}. Hence we omit the detail.
\qed
\end{proof}

\section{Microlocal index formula for local contributions}\label{sec:5}

In this section, using the Lefschetz cycle $LC(F,\Phi)$ introduced in Section \ref{sec:4}, we prove an index theorem which expresses local contributions of $(F,\Phi)$ as intersection numbers of the images of continuous sections of $\Lb \simeq T^*M$ and $LC(F,\Phi)$. Here we do not assume that the fixed point set $M$ of $\phi \colon X \longrightarrow X$ is smooth. However we assume the condition:
\begin{equation}
1 \notin \Ev_x \hspace{5mm} \text{for any $x \in M_{\reg}$}.
\end{equation}
Also in this more general setting, we can define the Lefschetz bundle $\Lb \simeq T^*M_{\reg}$ over $M_{\reg}$ and construct the Lefschetz cycle $LC(F,\Phi)$ in $\Lb$ by using the methods in Section \ref{sec:4}. Let $M = \bigsqcup_{i \in I}M_i$ be the decomposition of $M$ into connected components. Denote $(M_i)_{\reg}$ simply by $N_i$ and set $\Lb_i :=N_i \utimes{M_{\reg}}\Lb$. Then we get a decomposition $\Lb = \bigsqcup_{i \in I} \Lb_i \simeq \bigsqcup_{i \in I}T^*N_i$ of $\Lb$. By the direct sum decomposition
\begin{equation}
H^0_{\SS(F) \cap \Lb}(\Lb; \pi_{M_{\reg}}^{-1}\omega_{M_{\reg}}) \simeq \bigoplus_{i \in I}H^0_{\SS(F) \cap \Lb_i}(\Lb_i;\pi_{N_i}^{-1}\omega_{N_i}),
\end{equation}
we obtain a decomposition
\begin{equation}
LC(F,\Phi) =\dsum_{i \in I} LC(F,\Phi)_{M_i}
\end{equation}
of $LC(F,\Phi)$, where $LC(F,\Phi)_{M_i} \in H^0_{\SS(F) \cap \Lb_i}(\Lb_i;\pi_{N_i}^{-1}\omega_{N_i}) $. Now let us fix a fixed point component $M_i$ and assume that $\supp(F) \cap M_i$ is compact and contained in $N_i=(M_i)_{\reg}$. We shall show how the local contribution $c(F,\Phi)_{M_i}\in \Comp$ of $(F,\Phi)$ from $M_i$ can be expressed by $LC(F,\Phi)_{M_i}$. In order to state our results, for the sake of simplicity, we denote $N_i=(M_i)_{\reg}$, $\Lb_i$, $LC(F,\Phi)_{M_i}$, $c(F,\Phi)_{M_i}$ simply by $M$, $\Lb$, $LC(F,\Phi)$, $c(F,\Phi)$ respectively. Recall that to any continuous section $\sigma \colon M \longrightarrow \Lb \simeq T^*M$ of the vector bundle $\Lb$, we can associate a cycle $[\sigma] \in H^0_{\sigma(M)}(T^*M;\pi_M^!(\Comp_M))$ which is the image of $1 \in H^0(M;\Comp_M)$ by the isomorphism $H^0_{\sigma(M)}(T^*M;\pi_M^!\Comp_M) \simeq H^0(M;(\pi_M\circ \sigma)^!\Comp_M) \simeq H^0(M;\Comp_M)$ (see \cite[Definition 9.3.5]{K-S}). If $\sigma(M) \cap \supp (LC(F,\Phi))$ is compact, we can define the intersection number $\sharp ([\sigma] \cap LC(F,\Phi))$ of $[\sigma]$ and $LC(F,\Phi)$ to be the image of $[\sigma] \otimes LC(F,\Phi)$ by the chain of morphisms
\begin{eqnarray}
H^0_{\sigma(M)}(\Lb;\pi_M^!\Comp_M) \otimes H^0_{ \supp (LC(F,\Phi))}(\Lb;\pi_M^{-1}\omega_M)
&\longrightarrow & H^0_{\sigma(M) \cap \supp (LC(F,\Phi))} (\Lb;\omega_{\Lb}) \\
& \overset{\int_{\Lb}}{\longrightarrow} & \Comp.
\end{eqnarray}

\begin{thm}\label{thm:5-1}
Assume that $\supp (F) \cap M$ is compact. Then for any continuous section $\sigma \colon M \longrightarrow \Lb \simeq T^*M$ of $\Lb$, we have
\begin{equation}
c(F,\Phi) =\sharp ([\sigma] \cap LC(F,\Phi)).
\end{equation}
\end{thm}

\begin{proof}
Our proof is very similar to that of Kashiwara's microlocal index theorem (see \cite[Proposition 9.5.1]{K-S}). Set $S=\supp (F)$. Then the result follows from the commutative diagram \eqref{diag:5-8} below. By the commutativity of this diagram, the characteristic class $C(F,\Phi) \in H^0_{S\cap M}(M;\omega_M)$ and the Lefschetz cycle $LC(F,\Phi)\in H^0_{\SS(F) \cap \Lb}(\Lb; \pi_M^{-1}\omega_M)$ are sent to the same element in $H^0_{\pi_M^{-1}(S\cap M)} (\Lb; \pi_M^{-1}\omega_M)$ by the above morphisms {\bf A} and {\bf B}. Hence the proof proceeds just as in the way as that of \cite[Proposition 9.5.1]{K-S}. 
{\small \begin{equation}\label{diag:5-8}
\xymatrix@C=5mm{
R{\rm Hom}(F,F) \ar@{-}[d]^{\wr}& & & \\
\RG_{\Delta_X \cap (S\times S)}(X \times X;F \boxtimes \D F) \ar@{-}[d]^{\wr} \ar[r]& \RG_{\Delta_X \cap (S\times S)}(X \times X;h_*\omega_X) \ar@{-}[d]^{\wr} \ar@{-}[r]^{\sim} & \RG_{S\cap M}(M;\omega_M) \ar[d]^{\wr}_{{\bf A}}\\
\RG_{\pi_X^{-1}S}(T^*X ;\mu_{\Delta_X}(F \boxtimes \D F)) \ar[r] & \RG_{\pi_X^{-1}S}(T^*X;\mu_{\Delta_X}(h_*\omega_X)) \ar@{-}[r]^{\sim} & \RG_{\pi_M^{-1}(S\cap M)}(\Lb ; \pi_M^{-1}\omega_M) \\
\RG_{\SS(F)}(T^*X;\mu_{\Delta_X}(F \boxtimes \D F)) \ar[u]_{\wr} \ar[r]& \RG_{\SS(F)}(T^*X;\mu_{\Delta_X}(h_*\omega_X)) \ar[u] \ar@{-}[r]^{\sim}& \RG_{\SS(F) \cap \Lb}(\Lb ;\pi_M^{-1}\omega_M). \ar[u]^{{\bf B}}}
\end{equation}}
\qed
\end{proof}

As an application of Theorem \ref{thm:5-1}, we shall give a useful formula which enables us to describe the Lefschetz cycle $LC(F,\Phi)$ explicitly in the special case where $\phi \colon X \longrightarrow X$ is the identity map of $X$ and $M=X$. For this purpose, until the end of this section, we shall consider the situation where $\phi ={\rm id}_X$, $M=X$ and $\Phi \colon F \longrightarrow F$ is an endomorphism of $F \in \Dc(X)$. In this case, $LC(F,\Phi)$ is a Lagrangian cycle in $T^*X$. Now for real analytic function $\varphi \colon Y \longrightarrow I$ on a real analytic manifold $Y$ ($I$ is an open interval in $\Real$) we define a section $\sigma_{\varphi} \colon Y \longrightarrow T^*Y$ of $T^*Y$ by $\sigma_{\varphi}(y):=(y ;d\varphi(y))\ (y \in Y)$ and set 
\begin{equation}
\Lambda_{\varphi} :=\sigma_{\varphi}(Y)=\{ (y ;d\varphi(y)) \ |\ y \in Y \}. 
\end{equation}
Note that $\Lambda_{\varphi}$ is a Lagrangian submanifold of $T^*Y$ Then we have the following analogue of \cite[Theorem 9.5.3]{K-S}.

\begin{thm}\label{thm:5-2}
Let $Y$ be a real analytic manifold, $G$ an object of $\Dc(Y)$ and $\Psi \colon G \longrightarrow G$ an endomorphism of $G$. For a real analytic function $\varphi \colon Y \longtwoheadrightarrow I$, assume that the following conditions are satisfied.
\begin{enumerate}
\item $\supp (G) \cap \{ y \in Y\ |\ \varphi(y) \leq t\}$ is compact for any $t \in I$.
\item $\SS(G) \cap \Lambda_{\varphi}$ is compact.
\end{enumerate}
Then the global trace
\begin{equation}
\tr(G,\Psi)=\dsum_{j\in \Z}(-1)^j \tr \{ H^j(Y;G) \overset{\Psi}{\longrightarrow}H^j(Y;G)\}
\end{equation}
of $(G,\Psi)$ is equal to $\sharp ([\sigma_{\varphi}] \cap LC(G,\Psi))$.
\end{thm}

\begin{proof}
Since the fixed point set of $\phi =\id_Y$ is $Y$ itself, $LC(G,\Psi)$ is a Lagrangian cycle in $T^*Y$. Moreover, since any open subset of $Y$ is invariant by $\phi=\id_Y$, we can freely use the microlocal Morse lemma (\cite[Corollary 5.4.19]{K-S}) to reduce the computation of the global trace $\tr(G,\Psi)$ on $Y$ to that of
\begin{equation}
\dsum_{j\in \Z}(-1)^j \tr \{ H^j(\Omega_t;G) \overset{\Psi|_{\Omega_t}}{\longrightarrow}H^j(\Omega_t;G)\}
\end{equation}
for sufficiently large $t>0$ in $I$, where we set $\Omega_t:=\{y \in Y \ |\ \varphi(y)<t\}$. Then the proof proceeds essentially in the same way as that of \cite[Theorem 9.5.3]{K-S}.
\qed
\end{proof}

\begin{thm}\label{thm:5-3}
Let $X$, $F \in \Dc(X)$ and $\Phi \colon F \longrightarrow F$ be as above. For a real analytic function $\varphi \colon X \longrightarrow \Real$ and a point $x_0 \in X$, assume the condition
\begin{equation}
\Lambda_{\varphi} \cap \SS(F) \subset \{(x_0;d\varphi(x_0))\}.
\end{equation}
Then the intersection number $\sharp ([\sigma_{\varphi}] \cap LC(F,\Phi))$ (at the point $(x_0;d\varphi(x_0))\in T^*X$) is equal to
\begin{equation}
\dsum_{j\in \Z}(-1)^j \tr \{ H^j_{\{\varphi \geq \varphi(x_0)\}}(F)_{x_0} \overset{\Phi}{\longrightarrow}H^j_{\{\varphi \geq \varphi(x_0)\}}(F)_{x_0} \}.
\end{equation}
\end{thm}

\begin{proof}
The proof is very similar to that of \cite[Theorem 9.5.6]{K-S}. For a sufficiently small open ball $B(x_0,\e)=\{ x \in X \ |\ |x-x_0|<\e\}$ centered at $x_0$, set $F_0=\RG_{B(x_0,\e)}(F)\in \Dc(X)$. Then $\Phi$ induces a natural morphism $\Phi_0 \colon F_0 \longrightarrow F_0$ in $\Dc(X)$. Moreover by the proof of \cite[Theorem 9.5.6]{K-S}, we have
\begin{equation}\label{eq:5-13}
\Lambda_{\varphi} \cap \SS(F_0 ) \subset \pi_X^{-1}(\Omega_{-t}) \sqcup \{(x_0;d\varphi(x_0))\}
\end{equation}
for sufficiently small $t >0$, where we set $\Omega_k:=\{x \in X\ |\ \varphi(x) -\varphi(x_0) <k\}$ for $k \in \Real$. Then applying Theorem \ref{thm:5-2} to the case where $I=(-\infty, 0)$, $Y=\Omega_0$, $G=F_0|_{\Omega_0} \in \Dc(Y)$ and $\Psi=\Phi_0|_{\Omega_0} \colon G \longrightarrow G$, we obtain
\begin{eqnarray}
\lefteqn{\sharp ([\sigma_{\varphi}] \cap LC(F_0 ,\Phi_0 ) \cap \pi_X^{-1}(\Omega_0)) } \nonumber \\
&=& \dsum_{j \in \Z} (-1)^j \tr \{H^j(B(x_0,\e) \cap \Omega_0;F) \overset{\Phi}{\longrightarrow} H^j(B(x_0,\e) \cap \Omega_0;F) \}.\label{eq:5-14}
\end{eqnarray}
On the other hand, since $\supp (F_0 )$ is compact in $X$, by Theorem \ref{thm:5-1} we have
\begin{equation}\label{eq:5-15}
\sharp ([\sigma_{\varphi}] \cap LC(F_0 ,\Phi_0 ))= \dsum_{j \in \Z} (-1)^j\tr \{H^j(B(x_0,\e) ;F) \overset{\Phi}{\longrightarrow} H^j(B(x_0,\e);F) \}.
\end{equation}
Comparing \eqref{eq:5-14} with \eqref{eq:5-15} in view of \eqref{eq:5-13}, we see that the intersection number of $[\sigma_{\varphi}]$ and $LC(F_0 ,\Phi_0 )$ at $(x_0;d\varphi(x_0))$ is equal to
\begin{equation}\label{eq:5-17}
\dsum_{j \in \Z} (-1)^j \tr\{ H^j_{\{ \varphi \geq \varphi(x_0)\}} (F)_{x_0} \overset{\Phi}{\longrightarrow} H^j_{\{ \varphi \geq \varphi(x_0) \}}(F)_{x_0} \}.
\end{equation}
Since $LC(F,\Phi)= LC(F_0 ,\Phi_0 )$ in an open neighborhood of $(x_0;d\varphi(x_0))$ in $T^*X$, this last intersection number $\sharp ([\sigma_{\varphi}]\cap LC(F_0, \Phi_0))$ ($=$\eqref{eq:5-17}) is equal to $\sharp ([\sigma_{\varphi}] \cap LC(F,\Phi))$. This completes the proof. \qed
\end{proof}

By Theorem \ref{thm:5-3}, we can explicitly describe the Lefschetz cycle $LC(F,\Phi) \in \varGamma(T^*X;\L_X)$ as follows. Let $X= \bigsqcup_{\alpha \in A}X_{\alpha}$ be a $\mu$-stratification of $X$ such that
\begin{equation}
\supp(LC(F,\Phi)) \subset \SS(F) \subset \bigsqcup_{\alpha \in A} T_{X_{\alpha}}^*X.
\end{equation}
Then $\Lambda := \bigsqcup_{\alpha \in A} T_{X_{\alpha}}^*X$ is a closed conic subanalytic Lagrangian subset of $T^*X$. Moreover there exists an open dense smooth subanalytic subset $\Lambda_0$ of $\Lambda$ whose decomposition $\Lambda_0=\bigsqcup_{i \in I}\Lambda_i$ into connected components satisfies the condition
\begin{equation}\label{eq:5-18}
``\text{For any $i\in I$, there exists $\alpha_i \in A$ such that $\Lambda_i \subset T_{X_{\alpha_i}}^*X$. }"
\end{equation}

\begin{dfn}\label{dfn:5-4}
For $i \in I$ and $\alpha_i\in A$ as above, we define a complex number $m_i \in \Comp$ by
\begin{equation}\label{eq:5-19}
m_i :=\dsum_{j \in \Z} (-1)^j \tr\{H^j_{\{\varphi \geq \varphi(x)\}}(F)_x \overset{\Phi}{\longrightarrow} H^j_{\{\varphi \geq \varphi(x)\}}(F)_x\},
\end{equation}
where the point $x \in \pi_X(\Lambda_i) \subset X_{\alpha_i}$ and the $\Real$-valued real analytic function $\varphi \colon X \longrightarrow \Real$ (defined in an open neighborhood of $x$ in $X$) are defined as follows. Take a point $p \in \Lambda_i$ and set $x=\pi_X(p) \in X_{\alpha_i}$. Then $\varphi \colon X \longrightarrow \Real$ is a real analytic function which satisfies the following conditions:
\begin{enumerate}
\item $p=(x;d\varphi(x)) \in \Lambda_i$.
\item The Hessian ${\rm Hess} (\varphi|_{X_{\alpha_i}})$ of $\varphi|_{X_{\alpha_i}}$ is positive definite.
\end{enumerate}
\end{dfn}

\begin{cor}\label{cor:5-5}
In the situation as above, for any $i \in I$ there exists an open neighborhood $U_i$ of $\Lambda_i$ in $T^*X$ such that
\begin{equation}
LC(F,\Phi)=m_i \cdot [T_{X_{\alpha_i}}^*X]
\end{equation}
in $U_i$.
\end{cor}

Now let us define a $\Comp$-valued constructible function $\varphi(F,\Phi)$ on $X$ by
\begin{equation}
\varphi(F,\Phi)(x):=\dsum_{j \in \Z} (-1)^j\tr\{ H^j(F)_x \overset{\Phi|_{\{x\}}}{\longrightarrow} H^j(F)_x\}
\end{equation}
for $x\in X$. We will show that the characteristic cycle $CC(\varphi(F,\Phi))$ of $\varphi(F,\Phi)$ (see Proposition \ref{prp:2-10}) is equal to the Lefschetz cycle $LC(F,\Phi)$. For this purpose, we need the following.

\begin{dfn}[\cite{M-T-1}]
Let $\varphi \colon X \longrightarrow \Z$ be a $\Z$-valued constructible function on $X$ and $U$ a relatively compact subanalytic open subset in $X$. We define the topological integral $\dint_U \varphi$ of $\varphi$ over $U$ by
\begin{equation}
\dint_U \varphi =\dsum_{\alpha \in \Z} c_{\alpha} \cdot \chi(\RG(U;\Comp_{X_{\alpha}})),
\end{equation}
where $\varphi =\sum_{\alpha \in A} c_{\alpha} \1_{X_{\alpha}}$ ($c_{\alpha} \in \Z$) is an expression of $\varphi$ with respect to a subanalytic stratification $X=\bigsqcup_{\alpha \in A} X_{\alpha}$ of $X$.
\end{dfn}

We can extend $\Comp$-linearly this integral $\dint_U \colon \CF(X) \longrightarrow \Z$ and obtain a $\Comp$-linear map
\begin{equation}
\dint_U \colon \CF(X)_{\Comp} \longrightarrow \Comp.
\end{equation}
On the other hand, since any relatively compact subanalytic open subset $U$ of $X$ is invariant by $\phi=\id_X$, the global trace on $U$
\begin{equation}
\tr(F|_U,\Phi|_U)=\dsum_{j \in \Z}(-1)^j\tr\{H^j(U;F) \overset{\Phi|_U}{\longrightarrow}H^j(U;F)\}
\end{equation}
is well-defined.

\begin{lem}\label{lem:5-7}
For any relatively compact subanalytic open subset $U$ of $X$, we have
\begin{equation}
\tr(F|_U,\Phi|_U) =\dint_U \varphi(F,\Phi).
\end{equation}
\end{lem}

The proof of this lemma being completely similar to that of \cite[Proposition 11.6]{G-M-1}, we omit the proof.

\begin{thm}\label{thm:5-8}
In the situation $\phi=\id_X$, $\Phi \colon F \longrightarrow F$ etc. as above, we have the equality
\begin{equation}
LC(F,\Phi)=CC(\varphi(F,\Phi))
\end{equation}
as Lagrangian cycles in $T^*X$.
\end{thm}

\begin{proof}
Let $X=\bigsqcup_{\alpha \in A}X_{\alpha}$ be a $\mu$-stratification of $X$ such that
\begin{equation}
\supp (LC(F,\Phi)), \hspace{3mm}\supp(CC(\varphi(F,\Phi))) \subset \Lambda = \bigsqcup_{\alpha \in A}T_{X_{\alpha}}^*X.
\end{equation}
Take an open dense smooth subanalytic subset $\Lambda_0 $ of $\Lambda$ whose decomposition $\Lambda= \bigsqcup_{i \in I}\Lambda_i$ into connected components satisfies the condition \eqref{eq:5-18}. Let us fix $\Lambda_i$ and $X_{\alpha_i}$ such that $\Lambda_i \subset T_{X_{\alpha_i}}^*X$. It is enough to show that $LC(F,\Phi)$ and $CC(\varphi(F,\Phi))$ coincide in an open neighborhood of $\Lambda_i$ in $T^*X$. By Corollary \ref{cor:5-5}, in an open neighborhood $U_i$ of $\Lambda_i$ in $T^*X$ we have
\begin{equation}
LC(F,\Phi)=m_i \cdot [T_{X_{\alpha_i}}^*X],
\end{equation}
where $m_i \in \Comp$ is defined by \eqref{eq:5-19} for $p \in \Lambda_i$, $x=\pi_X(p)\in X_{\alpha_i}$, $\varphi \colon X \longrightarrow \Real$ as in Definition \ref{dfn:5-4}. Let $U$ be a sufficiently small open ball in $X$ centered at $x \in X_{\alpha_i}$. Set $V:=U \cap \{ \varphi <\varphi(x)\}$. Then we have
\begin{eqnarray}
m_i
&=& \dsum_{j \in \Z}(-1)^j\tr\{H^j_{\{\varphi \geq \varphi(x)\}}(U;F) \overset{\Phi}{\longrightarrow} H^j_{\{\varphi \geq \varphi(x)\}}(U;F) \} \\
&=& \tr(F|_U,\Phi|_U)-\tr(F|_V,\Phi|_V) \\
&=& \dint_U \varphi(F,\Phi) -\dint_V\varphi(F,\Phi).
\end{eqnarray}
This last number coincides with the coefficient of $[T_{X_{\alpha_i}}^*X]|_{U_i}$ in $CC(\varphi(F,\Phi))|_{U_i}$. This completes the proof. \qed
\end{proof}

\section{Explicit description of Lefschetz cycles}\label{sec:6}

In this section, we explicitly describe the Lefschetz cycle $LC(F,\Phi)$ introduced in Section \ref{sec:4} in many cases. Let $M$ be a possibly singular fixed point component of $\phi \colon X \longrightarrow X$. Throughout this section, we assume the condition
\begin{equation}
``\text{$1\notin \Ev_x$ for any $x \in \supp(F) \cap M_{\reg}$.}"
\end{equation}
Then there exists an open neighborhood $U$ of $\supp(F) \cap M_{\reg}$ in $M_{\reg}$ such that $\Gamma_{\phi}$ intersects with $\Delta_X$ cleanly along $U \subset M \subset \Gamma_{\phi} \cap \Delta_X$. Namely, there exists a Lefschetz bundle $\Lb=U \utimes{M}\{T^*_{\Gamma_{\phi}}(X\times X) \cap T^*_{\Delta_X}(X\times X)\}$ over $U$ which is isomorphic to $T^*U$. As in the same way as in Section \ref{sec:4}, we can define a Lagrangian cycle in $\Lb$ associated with $(F,\Phi)$. We still denote it by $LC(F,\Phi)$ and want to describe it explicitly. Replacing $X$, $M$ etc. by $X\setminus (M \setminus U)$, $U$ etc. respectively, we may assume that $M$ is smooth and $1 \notin \Ev_x$ for any $x \in M$ from the first. In this situation, the fixed point set of $\phi^{\prime} \colon T_MX \longrightarrow T_MX$ is the zero-section $M$. Let $\Gamma_{\phi^{\prime}}=\{(\phi^{\prime}(p),p) \ |\ p \in T_MX\}\subset T_MX \times T_MX$ be the graph of $\phi^{\prime}$ and $\Delta_{T_MX}\simeq T_MX$ the diagonal subset of $T_MX \times T_MX$. Then
\begin{equation}
\Lb^{\prime}:=T^*_{\Gamma_{\phi^{\prime}}}(T_MX\times T_MX) \cap T^*_{\Delta_{T_MX}}(T_MX \times T_MX)
\end{equation}
is a vector bundle over the zero-section $M \simeq \Gamma_{\phi^{\prime}} \cap \Delta_{T_MX}$ of $T_MX$. Since $\Lb^{\prime}$ is also isomorphic to $T^*M$ by our assumptions, we shall identify it with the original Lefschetz bundle $\Lb=T^*_{\Gamma_{\phi}}(X \times X) \cap T^*_{\Delta_X}(X\times X)$. Now consider the natural morphism
\begin{equation}
\Phi^{\prime} \colon (\phi^{\prime})^{-1}\nu_M(F) \longrightarrow \nu_M(F)
\end{equation}
induced by $\Phi \colon \phi^{-1}F \longrightarrow F$. Then from the pair $(\nu_M(F), \Phi^{\prime})$, we can construct the Lefschetz cycle $LC(\nu_M(F),\Phi^{\prime})$ in $\Lb' \simeq \Lb$.

\begin{prp}\label{prp:6-1}
In $\Lb \simeq \Lb^{\prime}$, we have
\begin{equation}
LC(F,\Phi) =LC(\nu_M(F),\Phi^{\prime}).
\end{equation}
\end{prp}

\begin{proof}
The proof is similar to those of \cite[Proposition 9.6.11]{K-S} and Proposition \ref{prp:3-3}. Indeed, the proof follows from the commutativity of Diagram 6.a, which is a microlocal version of Diagram \eqref{diag:3-6}. Here we denote $T_MX$, $\SS(F)$ and $C_{T_M^*X}(\SS(F))$ by $\Vb$, $S$ and $S^{\prime}$ respectively. Note that we have natural isomorphisms
\begin{equation}
T^*(T_MX)  \simeq  T^*(T^*_MX)  \simeq T_{T_M^*X}(T^*X)
\end{equation}
(see \cite[(6.2.3)]{K-S} and \eqref{eq:new} below) and the normal cone $S^{\prime}=C_{T_M^*X}(\SS(F))$ can be considered as a subset of $T^*(T_MX)=T^*\Vb$. We also used a conic isotropic subset $S^{\prime \prime}=(S \cap \Lb) \cup (S^{\prime} \cap \Lb^{\prime})$ of $\Lb \simeq \Lb^{\prime} \simeq T^*M$ and the morphism $\tl{h} \colon T_MX \longrightarrow T_MX \times T_MX$ is defined by $\tl{h}= (\phi^{\prime},\id_{T_MX})$. Moreover we used the natural isomorphism $\D \nu_M(F) \simeq \nu_M(\D F)$ to obtain Diagram 6.a. Let us explain the construction of the morphism ${\bf A}$ in Diagram 6.a. First consider the commutative diagram:
\begin{equation}
\xymatrix@R=3mm@C=5mm{
T_{M\times M}(X\times X) \ar@{^{(}->}[rr]^{s_1} & & \tl{(X\times X)_{M\times M}} & & \Omega_{X\times X} \ar@{_{(}->}[ll]_{j_1} \ar@{->>}[rr]^{\tl{p_1}}& & X\times X \\
 & \Box & & \Box & & \Box & \\
T_MX \ar@{^{(}->}[rr]^{s} \ar@{^{(}->}[uu]^{\delta_{T_MX}} & & \tl{X_M} \ar@{^{(}->}[uu]^{\tl{\delta^{\prime}}} & & \Omega_X \ar@{_{(}->}[ll]_{j} \ar@{->>}[rr]^{\tl{p}} \ar@{^{(}->}[uu]^{\tl{\delta}} & & X \ar@{^{(}->}[uu]^{\delta_X}}
\end{equation}
which already appeared in the proof of Proposition \ref{prp:3-3}. Denote the image of $\tl{\delta^{\prime}}$ (resp. $\tl{\delta}$) by $\Delta_{\tl{X_M}}$ (resp. $\Delta_{\Omega_X}$). Then we see that the following morphisms are isomorphisms. 
\begin{eqnarray}
\t{\tl{p_1}} &\colon &\Delta_{\Omega_X} \utimes{\Delta_X} T_{\Delta_X}^*(X\times X) \longrightarrow T_{\Delta_{\Omega_X}}^*\Omega_{X\times X}, \\
\t{j_1} &\colon &\Delta_{\Omega_X} \utimes{\Delta_{\tl{X_M}}}T_{\Delta_{\tl{X_M}}}^*(\tl{(X\times X)_{M\times M}}) \longrightarrow T_{\Delta_{\Omega_X}}^*\Omega_{X\times X}, \\
\t{s_1} &\colon &\Delta_{T_MX} \utimes{\Delta_{\tl{X_M}}}T_{\Delta_{\tl{X_M}}}^*(\tl{(X\times X)_{M\times M}}) \longrightarrow T_{\Delta_{T_MX}}^*
(T_{M\times M}(X\times X)). 
\end{eqnarray}

\newpage
\begin{center}{
\rotatebox[origin=c]{90}{
$
\xymatrix@R=10mm@C=10mm{
R{\rm Hom}(F,F) \ar[rr] & & R{\rm Hom}(\nu_M(F),\nu_M(F)) \\
\RG_{S}(T^*X;\mu_{\Delta_X}(F \boxtimes \D F)) \ar@{-}[u]^{\wr} \ar[r]^{{\bf A}} \ar[dd] & \RG_{S^{\prime}}(T^*\Vb;\mu_{\Delta_{\Vb}}(\nu_{M\times M}(F \boxtimes \D F))) \ar[d]& \RG_{S^{\prime}}(T^*\Vb;\mu_{\Delta_{\Vb}}(\nu_M(F) \boxtimes \D \nu_M(F))) \ar[l] \ar[d] \ar@{-}[u]^{\wr}\\
 & \RG_{S^{\prime}}(T^*\Vb;\mu_{\Delta_{\Vb}}(\tl{h}_*\tl{h}^{-1}\nu_{M\times M}(F \boxtimes \D F))) \ar[d]& \RG_{S^{\prime}}(T^*\Vb;\mu_{\Delta_{\Vb}}(\tl{h}_*(\phi^{\prime -1}\nu_M(F) \otimes \D \nu_M(F)))) \ar[l] \ar[d] \\
\RG_{S}(T^*X; \mu_{\Delta_X}(h_*(\phi^{-1}F \otimes \D F))  \ar[d]^{\Phi} \ar[r]^{{\bf B}} & \RG_{S^{\prime}}(T^*\Vb;\mu_{\Delta_{\Vb}}(\tl{h}_* \nu_M(\phi^{-1}F \otimes \D F))) \ar[d]^{\Phi} & \RG_{S^{\prime}}(T^*\Vb;\mu_{\Delta_{\Vb}}(\tl{h}_* (\nu_M(\phi^{-1}F) \otimes \D\nu_M(F)))) \ar[l] \ar[d]^{\Phi} \\
\RG_{S}(T^*X;\mu_{\Delta_X}(h_*(F\otimes \D F))) \ar[d] \ar[r] & \RG_{S^{\prime}}(T^*\Vb;\mu_{\Delta_{\Vb}}(\tl{h}_*\nu_M(F\otimes \D F)) \ar[d] & \RG_{S^{\prime}}(T^*\Vb;\mu_{\Delta_{\Vb}}(\tl{h}_*(\nu_M(F) \otimes \D \nu_M(F)))) \ar[l] \ar[d] \\
\RG_{S}(T^*X;\mu_{\Delta_X}(h_*\omega_X)) \ar[r] \ar[d] & \RG_{S^{\prime}}(T^*\Vb;\mu_{\Delta_{\Vb}}(\tl{h}_*\nu_M(\omega_X))) \ar[d] &  \RG_{S^{\prime}}(T^*\Vb;\mu_{\Delta_{\Vb}}(\tl{h}_*\omega_{\Vb})) \ar@{-}[l]^{\sim} \ar[d] \\
\RG_{S^{\prime \prime}}(\Lb; \pi_M^{-1}\omega_M) \ar@{=}[r]& \RG_{S^{\prime \prime}}(\Lb;\pi_M^{-1}\omega_M) \ar@{=}[r] & \RG_{S^{\prime \prime}}(\Lb;\pi_M^{-1}\omega_M).}
$}}

\bigskip
Diagram 6.a
\end{center}

\newpage
Now let us set 
\begin{eqnarray}
S_1&:=&\t{\tl{p_1}}(\Delta_{\Omega_X}\utimes{\Delta_X}S), \\
S_2&:=&\overline{\t{j_1}^{-1}(S_1)}, \\
S_3&:=&S_2 \cap T_{\Delta_{T_MX}}^*(T_{M\times M}(X\times X)).
\end{eqnarray}
Then we have the following morphisms
\begin{eqnarray}
\lefteqn{\RG_{\SS(F)} (T_{\Delta_X}^*(X\times X); \mu_{\Delta_X}(F \boxtimes \D F))} \nonumber \\
&\longrightarrow & \RG_{\SS(F)}(T_{\Delta_X}^*(X\times X); \mu_{\Delta_X}(R\tl{p_1}_*\tl{p_1}^{-1}(F \boxtimes \D F))) \\
&\longrightarrow & \RG_{S_1}(T_{\Delta_{\Omega_X}}^*\Omega_{X\times X}; \mu_{\Delta_{\Omega_X}}(\tl{p_1}^{-1}(F \boxtimes \D F))) \\
\label{eq:6-11}&\simot& \RG_{S_2}(T_{\Delta_{\tl{X_M}}}^*(\tl{(X\times X)_{M\times M}});\mu_{\Delta_{\tl{X_M}}}(Rj_{1*}\tl{p_1}^{-1}(F \boxtimes \D F)))\\
&\longrightarrow & \RG_{S_2}(T_{\Delta_{\tl{X_M}}}^*(\tl{(X\times X)_{M\times M}});\mu_{\Delta_{\tl{X_M}}} (s_{1*}s_1^{-1}Rj_{1*}\tl{p_1}^{-1}(F \boxtimes \D F))) \\
&\longrightarrow & \RG_{S_3}(T_{\Delta_{T_MX}}^*(T_{M\times M}(X\times X));\mu_{\Delta_{T_MX}}(s_1^{-1}Rj_{1*}\tl{p_1}^{-1}(F \boxtimes \D F)))\\
& =& \RG_{S_3}(T^* \Vb ;\mu_{\Delta_{\Vb}}(\nu_{M \times M}((F \boxtimes \D F)))),
\end{eqnarray}
where we used \cite[Theorem 4.3.2 and Proposition 3.3.9]{K-S} (see also the arguments in \cite[page 192-193]{K-S}) to prove that the morphism \eqref{eq:6-11} is an isomorphism. Let us show that $S_3$ is equal to $S^{\prime}$. Let $(x^{\prime},x^{\prime \prime})$ be a local coordinate system of $X$ such that $M=\{ x^{\prime}=0\}$ and $(x^{\prime},x^{\prime \prime}; \xi^{\prime},\xi^{\prime \prime})$ the associated coordinates of $T^*X$. Then by the Hamiltonian isomorphism etc., we can naturally identify $T^*(T_MX) \simeq T_{\Delta_{T_MX}}^* (T_{M\times M}(X\times X))$ with $T_{T_M^*X}(T^*X)$ as follows (see \cite[(6.2.3)]{K-S}). 
\begin{equation}\label{eq:new} 
\begin{array}{ccccc}
T^*(T_MX) & \simeq & T^*(T^*_MX) & \simeq & T_{T_M^*X}(T^*X). \\
\inun & & \inun & & \inun \\
(x^{\prime},x^{\prime \prime};\xi^{\prime}, \xi^{\prime \prime}) & \longleftrightarrow & (\xi^{\prime},x^{\prime \prime};-x^{\prime},\xi^{\prime \prime}) & \longleftrightarrow & (x^{\prime},x^{\prime \prime};\xi^{\prime}, \xi^{\prime \prime})
\end{array}
\end{equation}
Under this identification, we can prove that $S_3  \subset T^*(T_MX) \simeq T_{\Delta_{T_MX}}^*(T_{M\times M}(X\times X)) $ is equal to the normal cone $S^{\prime}= C_{T^*_MX}(\SS(F)) \subset T_{T^*_MX}(T^*X)$ as follows. In the associated local coordinates $(x^{\prime}, x^{\prime \prime},t;\xi^{\prime}, \xi^{\prime \prime})$ ($t>0$) of $\Delta_{\Omega_X} \utimes{\Delta_{\tl{X_M}}} T_{\Delta_{\tl{X_M}}}^*(\tl{(X\times X)_{M\times M}}) \ (\simeq \Delta_{\Omega_X} \utimes{\Delta_X} T_{\Delta_X}^*(X\times X) \simeq \Omega_X \utimes{X}T^*X)$, its subset $\t{j_1}^{-1}\t{\tl{p_1}}(\Delta_{\Omega_X} \utimes{\Delta_X}S)$ is expressed by
\begin{equation}
\{ (x^{\prime}, x^{\prime \prime},t;\xi^{\prime}, \xi^{\prime \prime}) \in \Delta_{\Omega_X} \utimes{\Delta_{\tl{X_M}}}T_{\Delta_{\tl{X_M}}}^*(\tl{(X\times X)_{M\times M}}) \ | \ (tx^{\prime}, x^{\prime \prime}; t^{-1}\xi^{\prime}, \xi^{\prime \prime}) \in \SS(F) \}.
\end{equation}
Hence we have 
\begin{eqnarray}
\lefteqn{(x^{\prime},x^{\prime \prime};\xi^{\prime},\xi^{\prime \prime}) \in S_3=S_2 \cap T_{\Delta_{T_MX}}^*(T_{M\times M}(X\times X))} \nonumber \\
&\Longleftrightarrow & \exists (x_n^{\prime}, x_n^{\prime \prime},t_n ;\xi_n^{\prime}, \xi_n^{\prime \prime}) \in \Delta_{\Omega_X} \utimes{\Delta_{\tl{X_M}}}T_{\Delta_{\tl{X_M}}}^*(\tl{(X\times X)_{M\times M}}) \nonumber \\
& & \hspace{10mm} \text{s.t.\ } \begin{cases} (x_n^{\prime}, x_n^{\prime \prime},t_n ;\xi_n^{\prime}, \xi_n^{\prime \prime}) \overset{n \to \infty}{\longrightarrow} (x^{\prime},x^{\prime \prime},0;\xi^{\prime},\xi^{\prime \prime}), \\ (t_nx_n^{\prime}, x_n^{\prime \prime}; t_n^{-1}\xi_n^{\prime}, \xi_n^{\prime \prime}) \in \SS(F) \end{cases}\\
&\Longleftrightarrow & \exists (x_n^{\prime}, x_n^{\prime \prime},t_n ;\xi_n^{\prime}, \xi_n^{\prime \prime}) \in \Delta_{\Omega_X} \utimes{\Delta_{\tl{X_M}}}T_{\Delta_{\tl{X_M}}}^*(\tl{(X\times X)_{M\times M}}) \nonumber \\
& & \hspace{10mm} \text{s.t.\ } \begin{cases} (x_n^{\prime}, x_n^{\prime \prime},t_n ;\xi_n^{\prime}, \xi_n^{\prime \prime}) \overset{n \to \infty}{\longrightarrow} (x^{\prime},x^{\prime \prime},0;\xi^{\prime},\xi^{\prime \prime}), \\ (t_nx_n^{\prime}, x_n^{\prime \prime}; \xi_n^{\prime}, t_n\xi_n^{\prime \prime}) \in \SS(F)
\end{cases}
\end{eqnarray}

\begin{eqnarray}
&\Longleftrightarrow & \exists ((\tl{x_n}^{\prime}, \tl{x_n}^{\prime \prime};\tl{\xi_n}^{\prime},\tl{\xi_n}^{\prime \prime}), c_n) \in \SS(F) \times \Real_{>0} \nonumber \\
& & \hspace{10mm}\text{s.t.\ }\begin{cases} (\tl{x_n}^{\prime}, \tl{x_n}^{\prime \prime};\tl{\xi_n}^{\prime},\tl{\xi_n}^{\prime \prime}) \overset{n \to \infty}{\longrightarrow} (0,x^{\prime \prime};\xi^{\prime}, 0), \\ (c_n\tl{x_n}^{\prime},c_n\tl{\xi_n}^{\prime \prime}) \overset{n \to \infty}{\longrightarrow} (x^{\prime}, \xi^{\prime \prime}) \end{cases} \\
&\Longleftrightarrow & (x^{\prime},x^{\prime \prime};\xi^{\prime},\xi^{\prime \prime}) \in S^{\prime}= C_{T^*_MX}(\SS(F)) \subset T_{T^*_MX}T^*X.
\end{eqnarray}
We thus obtained the morphism ${\bf A}$:
\begin{equation}
\RG_{S}(T^*X; \mu_{\Delta_X}(F \boxtimes \D F)) \longrightarrow \RG_{S^{\prime}}(T^*\Vb ;\mu_{\Delta_{\Vb}} (\nu_{M\times M}(F\boxtimes \D F))).
\end{equation}
We can construct also the morphism ${\bf B}$ in Diagram 6.a as follows.
\begin{eqnarray}
\lefteqn{\RG_S(T^*X; \mu_{\Delta_X}(h_*(\phi^{-1}F \otimes \D F)))} \nonumber \\
&\longrightarrow& \RG_{S^{\prime}}(T^*\Vb ; \mu_{\Delta_{\Vb}}(\nu_{M \times M}(h_*(\phi^{-1}F \otimes \D F)))) \\
&\longrightarrow& \RG_{S^{\prime}}(T^*\Vb ;\mu_{\Delta_{\Vb}}(\tl{h}_* \nu_M(\phi^{-1}F \otimes \D F))), 
\end{eqnarray}
where the first morphism is constructed in the same way as ${\bf A}$ and we used \cite[Proposition 4.2.4]{K-S} to construct the second morphism. This completes the proof. \qed
\end{proof}

Since we have
\begin{gather}
\phi|_M =\phi^{\prime}|_M=\id_M, \\
\Phi|_M=\Phi^{\prime}|_M \colon F|_M \longrightarrow F|_M
\end{gather}
($M$ is identified with the zero-section of $T_MX$), the Lefschetz cycle $LC(\nu_M(F)|_M,\Phi^{\prime}|_M)$ in $T^*M$ is the same as $LC(F|_M,\Phi|_M)$. In what follows, we shall identify $\Lb \simeq \Lb^{\prime}$ with $T^*M$ and compare $LC(F,\Phi)=LC(\nu_M(F),\Phi^{\prime})$ with $LC(F|_M,\Phi|_M)$.

\medskip
Since our result holds for any conic object on any vector bundle over $M$, let us consider the following general setting. Let $\pi \colon \Vb \longtwoheadrightarrow M$ be a real vector bundle over $M$ and $\psi \colon \Vb \longrightarrow \Vb$ an endomorphism of the vector bundle $\Vb$. Assume that the fixed point set of $\psi$ is the zero-section $M$ of $\Vb$. For each point $x \in M$ we define a finite subset $\Ev_x$ of $\Comp$ by
\begin{equation}
\Ev_x =\{ \text{ the eigenvalues of $\psi_x \colon \Vb_x \longrightarrow \Vb_x$}\} \subset \Comp
\end{equation}
as in the case of $\Vb=T_MX$ and $\psi=\phi^{\prime} \colon T_MX \longrightarrow T_MX$ (see Definition \ref{dfn:3-2}). Then the above assumption on the fixed point set of $\psi$ implies that $1 \notin \Ev_x$ for any $x \in M$. Suppose that we are given a conic $\Real$-constructible object $G \in \Dc(\Vb)$ on $\Vb$ and a morphism $\Psi \colon \psi^{-1}G \longrightarrow G$ in $\Dc(\Vb)$. From these data, we can construct the Lefschetz bundle $\Lb_0 \simeq T^*M$ associated with $\psi$ and the Lefschetz cycle $LC(G,\Psi)$ in it.

\begin{prp}\label{prp:6-2}
Let $x_0 \in M$ be a point of $M$ such that
\begin{equation}\label{eq:6-8}
\Ev_{x_0} \cap \Real_{\geq 1} =\emptyset.
\end{equation}
Then we have
\begin{equation}
LC(G,\Psi)=LC(G|_M,\Psi|_M)
\end{equation}
in an open neighborhood of $\pi_M^{-1}(x_0)$ in $\Lb_0 \simeq T^*M$.
\end{prp}

\begin{proof}
Take an open neighborhood $W$ of $x_0$ in $M$ such that $\Ev_x \cap \Real_{\geq 1} =\emptyset$ for any $x \in W$. Then there exists a closed ball
\begin{equation}
Z:=\overline{B(x_0,\e_0)} =\{x \in M \ |\ |x-x_0| \leq \e_0\} \hspace{5mm}(\e_0>0)
\end{equation}
in $W$ centered at $x_0$. Consider the conic object $G_{\pi^{-1}(Z)} \in \Dc(\Vb)$ and the morphism
\begin{equation}
\Psi_{\pi^{-1}(Z)} \colon \psi^{-1}(G_{\pi^{-1}(Z)}) \longrightarrow G_{\pi^{-1}(Z)}
\end{equation}
induced by $\Psi$. Since the construction of $LC(G,\Psi)$ and $LC(G|_M,\Psi|_M)$ is local and $x_0 \in {\rm Int}Z$, we may replace $(G,\Psi)$ by $(G_{\pi^{-1}(Z)},\Psi_{\pi^{-1}(Z)})$. By the homotopy invariance of $LC(G,\Psi)$ (see Proposition \ref{prp:4-9}), replacing $\psi$ by $\lambda \psi$ for $0<\lambda <1$ does not affect $LC(G,\Psi)$ nor $LC(G|_M,\Psi|_M)$. Hence by replacing $\psi$ by $\lambda \psi$ for sufficiently small $0<\lambda \ll1$, we may assume also that
\begin{equation}
\Ev_x \subset \{z \in \Comp\ |\ |z|<1\}
\end{equation}
for any $x \in \pi(\supp(G))$. Then there exists an open tubular neighborhood $D$ of the zero-section $M$ in $\Vb$ such that $\psi^{-1}(D)\supset D$, and we can construct a morphism
\begin{equation}
\RG_D(\Psi) \colon \psi^{-1}\RG_D(G) \longrightarrow \RG_D(G)
\end{equation}
induced by $\Psi \colon \psi^{-1}G \longrightarrow G$. Since $LC(\RG_D(G), \RG_D(\Psi))=LC(G,\Psi)$, we may replace the pair $(G,\Psi)$ by $(\RG_D(G),\RG_D(\Psi))$ and assume that $\supp(G)$ is compact. Let us take a $\mu$-stratification $\Vb= \bigsqcup_{\alpha \in A} \Vb_{\alpha}$ of $\Vb$ which satisfies the following three conditions.
\begin{enumerate}
\item There exists a subset $B \subset A$ such that the zero-section $M \subset \Vb$ of $\Vb$ is $\bigsqcup_{\beta \in B} \Vb_{\beta}$.
\item $\SS(G) \subset \bigsqcup_{\alpha \in A} T^*_{\Vb_{\alpha}}\Vb$ in $T^*\Vb$.
\item $\SS(G|_M) \subset \bigsqcup_{\beta \in B}T^*_{\Vb_{\beta}}M$ in $T^*M$.
\end{enumerate}
For $\beta \in B$, we shall denote $\Vb_{\beta} \subset M$ by $M_{\beta}$. Namely $M= \bigsqcup_{\beta \in B}M_{\beta}$ is a $\mu$-stratification of $M$. Set $\Lambda= \bigsqcup_{\beta \in B}T^*_{M_{\beta}}M \subset T^*M$. By the conditions above, we obtain
\begin{equation}
\supp(LC(G,\Psi)), \hspace{5mm}\supp(LC(G|_M,\Psi|_M)) \subset \Lambda.
\end{equation}
Therefore it suffices to show that $LC(G,\Psi)$ coincides with $LC(G|_M,\Psi|_M)$ on an open dense subset of $\Lambda$. Let $\Lambda_0$ be an open dense smooth subanalytic subset of $\Lambda$ whose decomposition $\Lambda_0= \bigsqcup_{i \in I}\Lambda_i$ into connected components satisfies the condition
\begin{equation}
``\text{For any $i \in I$, there exists $\beta_i \in B$ such that $\Lambda_i \subset T_{M_{\beta_i}}^*M$.}"
\end{equation}
Let us fix $\Lambda_i$ and $M_{\beta_i}$ as above and compare $LC(G,\Psi)$ with $LC(G|_M,\Psi|_M)$ on $\Lambda_i$. Take a point $p \in \Lambda_i$ and set $x =\pi_M(p) \in M_{\beta_i}$, where $\pi_M \colon T^*M \longrightarrow M$ is the projection. Let $\varphi \colon M \longrightarrow \Real$ be a real analytic function (defined in an open neighborhood of $x$) which satisfies that $p=(x;d\varphi(x)) \in \Lambda_i$, $\varphi(x)=0$ and the Hessian ${\rm Hess}(\varphi|_{M_{\beta_i}})$ is positive definite. Then by Corollary \ref{cor:5-5}, we have
\begin{equation}
LC(G|_M,\Psi|_M)=m_i \cdot [T_{M_{\beta_i}}^*M]
\end{equation}
in an open neighborhood of $\Lambda_i$ in $T^*M$, where $m_i \in \Comp$ is defined by
\begin{equation}
m_i:= \dsum_{j \in \Z} (-1)^j \tr\{ H^j_{\{ \varphi \geq 0\}}(B(x,\delta);G|_M) \overset{\Psi|_M}{\longrightarrow} H^j_{\{\varphi \geq 0\}}(B(x,\delta);G|_M)\}
\end{equation}
for sufficiently small $\delta >0$. Set $U:=B(x,\delta)$ and $V:=U \cap \{ \varphi<0\}$ in $M$. Then we have
\begin{equation}\label{eq:6-16}
m_i=\tr(\RG_U(G|_M),\RG_U(\Psi|_M))-\tr(\RG_V(G|_M),\RG_V(\Psi|_M)).
\end{equation}
Set also $\tl{U}:=\pi^{-1}(U)$, $\tl{V}:=\pi^{-1}(V)\subset \Vb$ and $\tl{\varphi}:=\varphi \circ \pi \colon \Vb \longrightarrow \Real$. Since $G$ is conic in an open neighborhood of the zero-section $M \subset \Vb$, we have
\begin{gather}
\RG_{\tl{U}}(G)|_M \simeq \RG_U(G|_M),\\
\RG_{\tl{V}}(G)|_M \simeq \RG_V(G|_M).
\end{gather}
Now let us set
\begin{gather}
\Lambda_{\varphi}:=\{(x;d\varphi(x))\ |\ x\in M\}\subset T^*M, \\
\Lambda_{\tl{\varphi}}:=\{(g;d\tl{\varphi}(g))\ |\ g\in \Vb\} \subset T^*\Vb.
\end{gather}
Then by Theorem \ref{thm:3-4}, it follows from our assumption \eqref{eq:6-8} for $x \in \supp(G) \cap M$ that
\begin{gather}
\tr(\RG_U(G|_M),\RG_U(\Psi|_M))=\tr(\RG_{\tl{U}}(G),\RG_{\tl{U}}(\Psi)), \\
\tr(\RG_V(G|_M),\RG_V(\Psi|_M))=\tr(\RG_{\tl{V}}(G),\RG_{\tl{V}}(\Psi)).
\end{gather}
Applying Theorem \ref{thm:5-1} to the pair $(\RG_{\tl{U}}(G),\RG_{\tl{U}}(\Psi))$, we obtain
\begin{equation}\label{eq:6-23}
\tr(\RG_U(G|_M),\RG_U(\Psi|_M))=\sharp ([\sigma_{\varphi}] \cap LC(\RG_{\tl{U}}(G), \RG_{\tl{U}}(\Psi))).
\end{equation}
Now by the condition (i) and the definition of $\Lambda$ we have
\begin{eqnarray}
\supp(LC(\RG_{\tl{U}}(G),\RG_{\tl{U}}(\Psi)))
&\subset & \SS(\RG_{\tl{U}}(G)) \cap \Lb_0 \label{eq:6-28}\\
&\subset & \{ \SS(G) \cup (SS(G) +T^*_{\partial \tl{U}}\Vb)\} \cap \Lb_0 \\
&\subset & \Lambda \cup (\Lambda +T^*_{\partial U}M)=:\Lambda^{\prime}.\label{eq:6-30}
\end{eqnarray}
Since $\Lambda^{\prime}$ is isotropic, by the microlocal Bertini-Sard theorem (\cite[Proposition 8.3.12]{K-S}) for $0<a\ll1$ we have
\begin{equation}\label{eq:6-31}
\Lambda^{\prime} \cap \Lambda_{\varphi} \cap \pi_M^{-1}(\{0 < | \varphi |< a \})=\emptyset.
\end{equation}
By the proof of \cite[Theorem 9.5.6]{K-S} (use \cite[(9.5.12) and (9.5.13)]{K-S}) and the estimate \eqref{eq:6-28}-\eqref{eq:6-30} and \eqref{eq:6-31}, shrinking $U=B(x,\delta)$ if necessary, we may assume from the first that
\begin{equation}
\Lambda_{\varphi} \cap \supp(LC(\RG_{\tl{U}}(G), \RG_{\tl{U}}(\Psi))) \subset \pi_M^{-1}(\{\varphi<-\e_0 \}) \sqcup \{p\}
\end{equation}
for sufficiently small $\e_0 >0$. 
Hence from \eqref{eq:6-23} we deduce
\begin{eqnarray}\label{eq:6-38}
\lefteqn{\tr(\RG_U(G|_M), \RG_U(\Psi|_M))} \nonumber \\
&=&\sharp \{\pi_M^{-1} (\{\varphi <-\e_0 \}) \cap [\sigma_{\varphi}] \cap LC(\RG_{\tl{U}}(G), \RG_{\tl{U}}(\Psi))\} +[\sigma_{\varphi}] \underset{p}{\cdot} LC(G,\Psi),
\end{eqnarray}
where $[\sigma_{\varphi}] \underset{p}{\cdot} LC(G,\Psi)$ is the local intersection number of $[\sigma_{\varphi}]$ and $LC(G,\Psi)$ at $p \in \Lambda_i$. The other term $\tr(\RG_V(G|_M),\RG_V(\Psi|_M))=\tr(\RG_{\tl{V}}(G),\RG_{\tl{V}}(\Psi))$ can be calculated as follows. For $\e>0$, set $V_{\e}:=V \cap \{ \varphi <-\e\}$ and $\tl{V_{\e}}:=\tl{V} \cap \{ \tl{\varphi}<-\e\}=\pi^{-1}(V_{\e})$.

\begin{lem}\label{lem:6-3}
For sufficiently small $\e >0$, we have
\begin{equation}
\tr(\RG_{\tl{V}}(G), \RG_{\tl{V}}(\Psi))=\tr(\RG_{\tl{V_{\e}}}(G),\RG_{\tl{V_{\e}}}(\Psi)).
\end{equation}
\end{lem}

\begin{proof}
Set $\Sigma:=\SS(\RG_{\tl{U}}(G)) \subset T^*\Vb$. Then by the microlocal Bertini-Sard theorem (\cite[Proposition 8.3.12]{K-S}) there exists $\e>0$ such that
\begin{equation}
\Sigma \cap \Lambda_{\tl{\varphi}} \cap \pi^{-1}(\{-\e <\tl{\varphi}<0\} )=\emptyset.
\end{equation}
Hence by \cite[Corollary 5.4.19]{K-S}, we obtain
\begin{equation}
\RG(\{\tl{\varphi}<0\};\RG_{\tl{U}}(G)) \simto \RG(\{\tl{\varphi}<-\e\};\RG_{\tl{U}}(G)).
\end{equation}
\qed
\end{proof}

Let us continue the proof of Proposition \ref{prp:6-2}. By Lemma \ref{lem:6-3} and Theorem \ref{thm:5-1}, we obtain
\begin{equation}\label{eq:6-27}
\tr(\RG_V(G|_M),\RG_V(\Psi|_M))=\sharp ([\sigma_{\varphi}] \cap LC(\RG_{\tl{V_{\e}}}(G),\RG_{\tl{V_{\e}}}(\Psi)))
\end{equation}
for sufficiently small $\e>0$. Moreover it follows from the condition (i) and the definition of $\Lambda$ that
\begin{eqnarray}
\supp(LC(\RG_{\tl{V_{\e}}}(G),\RG_{\tl{V_{\e}}}(\Psi)))
&\subset & \SS(\RG_{\{\tl{\varphi}<-\e\}}(\RG_{\tl{U}}(G))) \cap \Lb_0 \\
&\subset & \Lambda^{\prime} +\Real_{\leq 0}\Lambda_{\varphi}.
\end{eqnarray}
Comparing this last estimate with \eqref{eq:6-31}, we obtain
\begin{equation}
\Lambda_{\varphi} \cap \supp(LC(\RG_{\tl{V_{\e}}}(G),\RG_{\tl{V_{\e}}}(\Psi))) \subset \pi_M^{-1}(\{\varphi<-\e_0 \})
\end{equation}
for $0<\e \ll \e_0$. Since
\begin{equation}
LC(\RG_{\tl{V_{\e}}}(G), \RG_{\tl{V_{\e}}}(\Psi))=LC(\RG_{\tl{U}}(G),\RG_{\tl
{U}}(\Psi))
\end{equation}
on $\pi_M^{-1}(\{\varphi<-\e_0 \})$, from \eqref{eq:6-27} we obtain
\begin{equation}\label{eq:6-36}
\tr(\RG_V(G|_M),\RG_V(\Psi|_M))= \sharp \{\pi_M^{-1} (\{\varphi<-\e_0 \}) \cap [\sigma_{\varphi}] \cap LC(\RG_{\tl{U}}(G), \RG_{\tl{U}}(\Psi))\}. 
\end{equation}
Putting \eqref{eq:6-38} and \eqref{eq:6-36} into \eqref{eq:6-16}, we finally obtain
\begin{equation}
m_i =[\sigma_{\varphi}] \underset{p}{\cdot} LC(G,\Psi),
\end{equation}
which shows
\begin{equation}
LC(G,\Psi)=LC(G|_M, \Psi|_M)
\end{equation}
on $\Lambda_i$. This completes the proof.
\qed
\end{proof}

Combining Proposition \ref{prp:6-1} and \ref{prp:6-2} with Theorem \ref{thm:5-8}, we can obtain explicit descriptions of the Lefschetz cycle $LC(F,\Phi)$ as follows. Let $\varphi(F|_M, \Phi|_M)$ be a $\Comp$-valued constructible function on $M$ defined by 
\begin{equation}
\varphi(F|_M, \Phi|_M)(x) =\dsum_{j \in \Z} (-1)^j\tr\{ H^j(F)_x \overset{\Phi|_{\{x\}}}{\longrightarrow} H^j(F)_x\}
\end{equation}
for $x\in M$. 

\begin{thm}\label{thm:6-4}
Let $x_0 \in M$ be a point of $M$ such that
\begin{equation}\label{eq:6-42}
\Ev_{x_0} \cap \Real_{\geq 1}=\emptyset.
\end{equation}
Then we have
\begin{equation}
LC(F,\Phi)=LC(F|_M,\Phi|_M)=CC(\varphi(F|_M,\Phi|_M))
\end{equation}
in an open neighborhood of $\pi_M^{-1}(x_0)$ in $T^*M$. 
\end{thm}

In the complex case, we have the following stronger result.

\begin{thm}\label{thm:6-5}
In the situation as above, assume moreover that $X$ and $\phi \colon X \longrightarrow X$ are complex analytic and $F \in \Db_c(X)$ i.e. $F$ is $\Comp$-constructible. Then we have
\begin{equation}
LC(F,\Phi)=LC(F|_M,\Phi|_M)=CC(\varphi(F|_M, \Phi|_M))
\end{equation}
globally on $T^*M$.
\end{thm}

\begin{proof}
By Proposition \ref{prp:6-1}, we have only to prove
\begin{equation}\label{eq:6-43}
LC(\nu_M(F),\Phi^{\prime})=LC(F|_M,\Phi|_M).
\end{equation}
Since these cycles are considered as sections of the sheaf of $\L_M$ of Lagrangian cycles on $T^*M$, it suffices to prove \eqref{eq:6-43} locally. Namely, for each $x_0 \in M$ we have only to prove \eqref{eq:6-43} in an open neighborhood of $\pi_M^{-1}(x_0)$ in $\Lb \simeq T^*M$. This local statement can be proved along the same line as the proof of Proposition \ref{prp:6-2}. Since $\nu_M(F)$ admits the action of $\Comp^{\times}$ in the complex case, we may use the arguments in the proof of \cite[Corollary 9.6.16]{K-S} for this purpose. This completes the proof.
\qed
\end{proof}

By Theorem \ref{thm:6-5} above, we can drop the assumption of the smoothness of $M$ or $\supp(F) \cap M$ in Theorem \ref{thm:3-5} (we can also drop the assumption (iii) of Corollary \ref{cor:3-7}).

\begin{cor}\label{cor:6-6}
Let $X$, $\phi$ and $M$ be as above and $F_1 \overset{\alpha}{\longrightarrow} F_2 \overset{\beta}{\longrightarrow} F_3 \overset{\gamma}{\longrightarrow}+1$ a distinguished triangle in $\Dc(X)$. Assume that we are given a morphism of distinguished triangles
\begin{equation}
\xymatrix@R=8mm@C=10mm{
\phi^{-1}F_1 \ar[r]^{\phi^{-1}\alpha} \ar[d]^{\Phi_1} & \phi^{-1}F_2 \ar[r]^{\phi^{-1}\beta} \ar[d]^{\Phi_2} & \phi^{-1}F_3 \ar[r]^{\phi^{-1}\gamma} \ar[d]^{\Phi_3} & \phi^{-1}F_1[1] \ar[d]^{\Phi_1[1]}\\
F_1 \ar[r]^{\alpha} & F_2 \ar[r]^{\beta} &F_3 \ar[r]^{\gamma} & F_1[1] }
\end{equation}
in $\Dc(X)$. Then for any $x_0 \in M$ such that $\Ev_{x_0} \cap \Real_{\geq 1}=\emptyset$, we have
\begin{equation}
LC(F_2,\Phi_2)=LC(F_1,\Phi_1)+LC(F_3,\Phi_3)
\end{equation}
in an open neighborhood of $\pi_M^{-1}(x_0)$ in $T^*M$.
\end{cor}

\section{Another construction of Lefschetz cycles}\label{sec:7}

In this section, we shall introduce another construction of Lefschetz cycles which slightly differs from the previous one. Moreover we prove that the difference is expressed by the sign $\pm 1$ of the determinant of $\id-\phi_x^{\prime} \colon (T_MX)_x \longrightarrow (T_MX)_x$ for $x \in M$. Since (except Proposition \ref{prp:7-1} below) the results in this section will be used only in the proof of our inverse image theorem in Section \ref{sec:8}, the readers who do not require the inverse image theorem can skip this section.

\subsection{New construction of Lefschetz cycles}\label{sec:7-1}

In this subsection, we inherit the situation and notations in previous sections and consider the problem in an open neighborhood $U$ of a smooth point of $M$ for which the condition 
\begin{equation}
1 \notin \Ev_x \hspace{5mm} \text{for any $x \in M \cap U$}
\end{equation}
is satisfied. Then we can construct locally the Lefschetz bundle $\Lb$ over $M$. Before introducing another construction of Lefschetz cycles, first let us study the structure of the object 
\begin{equation}
G:=\nu_{\Delta_X}(h_*\omega_X) |_{M \utimes{X}TX} \in \Dc(M \utimes{X}TX). 
\end{equation}
Let $g \colon \Gamma_{\phi} \longhookrightarrow X \times X$ be the inclusion map of the graph of $\phi$. Then we obtain an injective map
\begin{equation}\label{eq:7-3}
g^{\prime} \colon T_MX \simeq T_M \Gamma_{\phi} \longhookrightarrow T_{\Delta_X}(X \times X) \simeq TX
\end{equation}
induced by $g$. Recall that in Section \ref{sec:4} we defined a subbundle $\tLb \subset M \utimes{X}TX$ to be the image of this map. Let $i_{\tLb} \colon \tLb \longhookrightarrow M \utimes{X}TX$ be the inclusion map.

\begin{prp}\label{prp:7-1}
In the situation as above, we have an isomorphism 
\begin{equation}
G \simeq (i_{\tLb})_* \omega_{\tLb}.
\end{equation}
\end{prp}

\begin{proof}
Consider the following standard commutative diagram for the normal deformation $\tl{(\Gamma_{\phi})_M}$ of $\Gamma_{\phi}$ along $M \simeq \Delta_M \subset \Gamma_{\phi}$:
\begin{equation}
\xymatrix@C=15mm{
T_M\Gamma_{\phi} \ar@{^{(}->}[r]^{s_2} \ar[d] & \tl{(\Gamma_{\phi})_M} \ar[d]_{p_2} & \Omega_{\Gamma}=\{t_2 >0\} \ar@{_{(}->}[l]_{j_2} \ar@{->>}[ld]^{\tl{p_2}} \\
M\simeq \Delta_M \ar@{^{(}->}[r] &\Gamma_{\phi},}
\end{equation}
where $t_2 \colon \tl{(\Gamma_{\phi})_M} \longrightarrow \Real$ is the deformation parameter. Then we have the following Cartesian diagrams
\begin{equation}
\xymatrix@R=3mm@C=5mm{
T_M\Gamma_{\phi} \ar@{^{(}->}[rr]^{s_2} \ar@{^{(}->}[dd]^{g^{\prime}} & & \tl{(\Gamma_{\phi})_M} \ar@{^{(}->}[dd]^{\tl{g^{\prime}}}& & \Omega_{\Gamma} \ar@{_{(}->}[ll]_{j_2} \ar@{^{(}->}[dd]^{\tl{g}} \ar@{->>}[rr]^{\tl{p_2}}& & \Gamma_{\phi} \ar@{^{(}->}[dd]^g \\
 & \Box & & \Box & & \Box & \\
T_{\Delta_X}(X\times X) \ar@{^{(}->}[rr]^{s_0} & & \tl{(X\times X)_{\Delta_X}}& & \Omega \ar@{_{(}->}[ll]_{\hspace{5mm}j_0} \ar@{->>}[rr]^{\tl{p_0}} & & X\times X }
\end{equation}
induced by $g \colon \Gamma_{\phi} \longhookrightarrow X \times X$. From this we obtain an isomorphism
\begin{equation}
(g^{\prime})_*\nu_M(\omega_{\Gamma_{\phi}}) \simeq \nu_{\Delta_X}(g_*\omega_{\Gamma_{\phi}}) \simeq \nu_{\Delta_X}(h_*\omega_X). 
\end{equation}
Since we have $\nu_M(\omega_{\Gamma_{\phi}}) \simeq \omega_{T_M\Gamma_{\phi}}$ and $T_M\Gamma_{\phi}$ is identified with $\tLb$ by $g^{\prime}$, the result follows. 
\qed
\end{proof}

From now on, we shall introduce another construction of Lefschetz cycles. Let $\iota \colon TM \longhookrightarrow M \utimes{X}TX$ be the natural injection and $\rho \colon M \utimes{X}T^*X \longtwoheadrightarrow T^*M \simeq \Lb$ its dual. Let $i_0 \colon M \longhookrightarrow TM$ be the zero-section embedding. Then we have an isomorphism
\begin{equation}\label{eq:7-13}
R\rho_!(G^{\wedge}) \simeq (\iota^!\Comp_{M \utimes{X}TX} \overset{L}{\otimes} \iota^{-1}G)^{\wedge}
\end{equation}
in $\Dc(T^*M)$ by \cite[Proposition 3.7.14]{K-S}, where $( \ \cdot \ )^{\wedge}$ stands for the Fourier-Sato transform. Note that we have $G^{\wedge} \simeq \mu_{\Delta_X} (h_*\omega_X)|_{M \utimes{X}T^*X}$ and $\rho$ is proper on $\supp (G^{\wedge}) =\Lb$. The structure of the right hand side of the isomorphism \eqref{eq:7-13} is given by the next lemma.

\begin{lem}\label{lem:7-2}
We have a natural isomorphism
\begin{equation}
(\iota^!\Comp_{M \utimes{X}TX} \overset{L}{\otimes} \iota^{-1}G)^{\wedge}\simeq \pi_M^{-1}\omega_M.
\end{equation}
\end{lem}

\begin{proof}
Since $\iota^{-1}(\supp(G)) =\iota^{-1}(\tLb)$ is the zero-section $M$ of $TM$, we have
\begin{equation}
\iota^{-1}G \simeq i_{0*}(G|_M) \simeq  i_{0*}\{ (h_*\omega_X)|_{\Delta_M}\} \simeq  i_{0*}(\omega_X|_M), 
\end{equation}
where 
$M \simeq \Delta_M$ is identified with the zero-section of $M \utimes{X}TX \simeq \Delta_M \utimes{\Delta_X}T_{\Delta_X}(X\times X)$. Hence we obtain 
\begin{equation}\label{eq:7-20}
\iota^!\Comp_{M \utimes{X}TX} \overset{L}{\otimes} \iota^{-1}G \simeq i_{0*}(\omega_{M/X} \overset{L}{\otimes} \omega_X|_M) \simeq  i_{0*} \omega_M
\end{equation}
and 
\begin{equation}
(\iota^!\Comp_{M \utimes{X}TX} \overset{L}{\otimes} \iota^{-1}G)^{\wedge} \simeq  (i_{0*}\omega_M)^{\wedge} \simeq  \pi_M^{-1}\omega_M.
\end{equation}
\qed
\end{proof}

By Lemma \ref{lem:7-2}, we have a chain of morphisms
\begin{eqnarray}
R{\rm Hom}_{\Comp_X}(F,F)
&\simeq & \RG_{\SS(F)} (T^*X;\mu_{\Delta_X}(F \boxtimes DF)) \\
&\longrightarrow & \RG_{\SS(F)}(T^*X;\mu_{\Delta_X}(h_*\omega_X)) \\
&\simeq & \RG_{\SS(F) \cap \Lb} (M \utimes{X}T^*X; G^{\wedge})\\
&\simeq & \RG_{\SS(F) \cap \Lb} (T^*M;R\rho_!(G^{\wedge}))
\label{eq:7-27}\\
&\simeq & \RG_{\SS(F) \cap \Lb} (T^*M;
\pi_M^{-1}\omega_M),\label{eq:7-27-2}
\end{eqnarray}
where for the construction of the second (resp. last) morphism we used the morphisms \eqref{eq:4-6}-\eqref{eq:4-11} in the construction of $LC(F,\Phi)$ of Section \ref{sec:4} (resp. Lemma \ref{lem:7-2} and the isomorphism \eqref{eq:7-13}). By taking the $0$-th cohomology groups of both sides, we obtain a morphism
\begin{equation}\label{eq:7-28}
{\rm Hom}_{\Db(X)}(F,F) \longrightarrow H^0_{\SS(F) \cap \Lb}(T^*M;\pi_M^{-1}\omega_M).
\end{equation}

\begin{dfn}\label{dfn:7-3}
We denote by $\overline{LC(F,\Phi)}$ the image of $\id_F \in {\rm Hom}_{\Db(X)}(F,F)$ in \\$H^0_{\SS(F) \cap \Lb}(T^*M;\pi_M^{-1}\omega_M)$ by the morphism \eqref{eq:7-28}.
\end{dfn}

\subsection{Relations between $LC(F,\Phi)$ and $\overline{LC(F,\Phi)}$}\label{sec:7-3}

In this subsection, we shall compare $LC(F,\Phi)$ with $\overline{LC(F,\Phi)}$. For this purpose, we first consider the isomorphism
\begin{equation}
R\rho_*(G^{\wedge}) \simeq (\iota^!G)^{\wedge}
\end{equation}
obtained by \cite[Proposition 3.7.14]{K-S}. The right hand side is calculated as follows. Since we have
\begin{eqnarray}
\iota^!G
&\simeq & i_{0*}(\RG_M(G)|_M) \label{eq:7-31}\\
&\simeq & i_{0*}(\RG_{\Delta_X}(\nu_{\Delta_X}(h_*\omega_X))|_M) \\
&\simeq & i_{0*}(\nu_{\Delta_X} ( \RG_{\Delta_X} (h_*\omega_X))|_{\Delta_M}) \\
&\simeq & i_{0*}\omega_M,\label{eq:7-33}
\end{eqnarray}
we obtain an isomorphism 
\begin{equation}
(\iota^!G)^{\wedge} \simeq  (i_{0*}\omega_M)^{\wedge} \simeq  \pi_M^{-1}\omega_M.
\end{equation}
By using the isomorphism 
\begin{equation}
R\rho_*(G^{\wedge}) \simeq (\iota^!G)^{\wedge} \simeq  \pi_M^{-1}\omega_M
\end{equation}
thus obtained to change the construction of the morphism from  \eqref{eq:7-27} to \eqref{eq:7-27-2}, we obtain also a morphism
\begin{equation}
{\rm Hom}_{\Db(X)}(F,F) \longrightarrow H^0_{\SS(F) \cap \Lb}(T^*M;\pi_M^{-1}\omega_M).
\end{equation}
Then the image of $\id_F$ by this morphism is $LC(F,\Phi)$. Since $\supp(G) =\tLb \subset M\utimes{X}TX$ and $\iota \colon TM \longhookrightarrow M \utimes{X}TX$ is non-characteristic for $G$, we obtain an isomorphism
\begin{equation}
\gamma \colon \iota^!\Comp_{M\utimes{X}TX} \overset{L}{\otimes}\iota^{-1}G \simto \iota^!G
\end{equation}
by \cite[Proposition 5.4.13]{K-S}. Note that we have the following commutative diagram.
\begin{equation}
\xymatrix@C=20mm{
R \rho_! (G^{\wedge}) \ar[r]^{\sim} \ar[d]_{\wr} & (\iota^!\Comp_{M\utimes{X}TX} \overset{L}{\otimes}\iota^{-1}G)^{\wedge} \ar[d]^{\gamma}_{\wr} \\
R \rho_* (G^{\wedge}) \ar[r]_{\sim} & (\iota^!G)^{\wedge}.}
\end{equation}
Recall that for the constructions of $\overline{LC(F,\Phi)}$ and $LC(F,\Phi)$ we used the isomorphisms 
\begin{gather}
\alpha \colon \iota^!\Comp_{M\utimes{X}TX} \overset{L}{\otimes}\iota^{-1}G \simto i_{0*}\omega_M, \\
\beta \colon \iota^!G \simto i_{0*}\omega_M
\end{gather}
obtained in \eqref{eq:7-20} and \eqref{eq:7-31}-\eqref{eq:7-33} respectively.  However the diagram
\begin{equation}
\xymatrix@C=20mm{
\iota^!\Comp_{M\utimes{X}TX} \overset{L}{\otimes}\iota^{-1}G \ar[d]^{\gamma}_{\wr} \ar[r]^{\hspace*{3mm}\alpha}_{\hspace*{3mm}\sim} & i_{0*}\omega_M \\
\iota^!G \ar[ur]^{\beta}_{\sim} &}
\end{equation}
is not commutative in general. Hence our new Lefschetz cycle $\overline{LC(F,\Phi)}$ may be different from the original one $LC(F,\Phi)$. In order to describe the difference, from now on we shall assume that $M$ is connected. Then the following definition makes sense.

\begin{dfn}\label{dfn:7-4}
Define a number $\sgn(\id-\phi^{\prime})\in \{\pm 1\}$ to be the sign of the determinant of
\begin{equation}
\id-\phi^{\prime}_x \colon (T_MX)_x \longrightarrow (T_MX)_x,
\end{equation}
where $x$ is a point in $M$.
\end{dfn}

\begin{prp}\label{prp:7-5}
The following diagram is commutative.
\begin{equation}
\xymatrix@C=20mm{
\iota^!\Comp_{M\utimes{X}TX} \overset{L}{\otimes}\iota^{-1}G \ar[r]^{\alpha}_{\sim} \ar[d]^{\gamma}_{\wr} & i_{0*}\omega_M \ar[d]^{\sgn(\id-\phi^{\prime}) \times }_{\wr} \\
\iota^!G \ar[r]^{\beta}_{\sim} &i_{0*}\omega_M.}
\end{equation}
\end{prp}

For the proof of Proposition \ref{prp:7-5}, we need some refined arguments on orientation sheaves. For this purpose, we first prepare two key lemmas (Lemma \ref{lem:7-6} and \ref{lem:7-7} below) concerning orientation sheaves. Now let $Y$ be an $n$-dimensional real analytic manifold and $N \subset Y$ a submanifold of codimension $d$. Then we have the canonical isomorphism
\begin{equation}
or_N \simeq \H_N^d(or_Y)|_N
\end{equation}
for orientation sheaves. Let us describe this isomorphism more explicitly in terms of differential forms with hyperfunction coefficients. 

First recall that for an open subset $U \subset N$ there exists an isomorphism
\begin{equation}
or_N(U) \simeq [H_c^{n-d}(U;\Comp_N)]^*.
\end{equation}
Therefore to an orientation $\sigma_N$ of $N$ we can associate a section $1_{\sigma_N} \in or_N(U)$ which corresponds to the linear map
\begin{equation}
\begin{array}{ccc}
H_c^{n-d}(U;\Comp_N) & \longrightarrow & \Comp, \\
\inun & & \inun \\
{[\omega]}  & \longmapsto & \dint_{U_{\sigma_N}}1\cdot \omega
\end{array}
\end{equation}
where $\omega$ is a $C^{\infty}$-differential $(n-d)$-form with compact support on $U$ and $\dint_{U_{\sigma_N}}$ stands for the integration over $U$ with respect to the orientation $\sigma_N$.

On the other hand, a section of the sheaf $\H_N^d(or_Y)|_N$ over $U \subset N$ can be explicitly expressed as follows. Let $\B_Y^{(i)}$ be the sheaf of differential $i$-forms on $Y$ with hyperfunction coefficients. Then we have a flabby resolution of $or_Y$:
\begin{equation}
0 \longrightarrow or_Y \longrightarrow or_Y \uotimes{\Comp}\B_Y^{(0)} \underset{d}{\longrightarrow} or_Y \uotimes{\Comp} \B_Y^{(1)} \underset{d}{\longrightarrow} or_Y \uotimes{\Comp} \B_Y^{(2)} \underset{d}{\longrightarrow} \cdots.
\end{equation}
By taking an open subset $\tl{U}$ in $Y$ such that $\tl{U} \cap N=U$, we obtain an isomorphism
\begin{eqnarray}
\lefteqn{\H_N^d(or_Y)|_N(U)} \nonumber \\
&\simeq& H^d[ 0 \longrightarrow (or_Y \uotimes{\Comp}\varGamma_N\B_Y^{(0)})(\tl{U}) \underset{d}{\longrightarrow} \cdots \underset{d}{\longrightarrow} (or_Y \uotimes{\Comp}\varGamma_N\B_Y^{(d)})(\tl{U}) \longrightarrow \cdots ].
\end{eqnarray}
In this way, a section of $\H_N^d(or_Y)|_N$ over $U$ can be represented by an element of $(or_Y \uotimes{\Comp} \varGamma_N\B_Y^{(d)})(\tl{U})$. Now let $y=(y_1,\cdots,y_n)$ be a local coordinate system of $Y$ such that $N=\{y_1=\cdots =y_d=0\}$. Set $y^{\prime}=(y_1,\cdots, y_d)$ and denote by $\delta (y^{\prime}) \in \varGamma_N\B_Y^{(0)}$ Dirac's delta function on $Y$ supported by $N$. Then we have the following. 

\begin{lem}\label{lem:7-6}
By the isomorphism $or_N(U) \simeq \H_N^d(or_Y)|_N(U)$, the section $1_{\sigma_N}\in or_N(U)$ corresponds to $[\omega ]\in \H_N^d(or_Y)|_N(U)$, where $\omega \in (or_Y \uotimes{\Comp} \varGamma_N\B_Y^{(d)})(\tl{U})$ is defined by
\begin{equation}
\omega =1_{dy_1 \wedge \cdots \wedge dy_d \wedge \sigma_N} \otimes \delta (y^{\prime}) dy_1 \wedge \cdots \wedge dy_d.
\end{equation}
\end{lem}

The proof of this lemma follows from the definition of orientation sheaves (the proof of the Poincar{\'e}-Verdier duality theorem). Since it seems that this lemma is well-known to specialists, we omit the proof here. Similarly we have

\begin{lem}\label{lem:7-7}
By the isomorphism $or_{N/Y}(U) \simeq \H_N^d(\Comp_Y)|_N(U)$, the section \\
$1_{\sigma_N} \otimes (1_{dy_1 \wedge \cdots \wedge dy_d \wedge \sigma_N})^{\otimes -1} \in or_{N/Y}(U)$ corresponds to $[\omega_0] \in \H_N^d(\Comp_Y)|_N(U)$, where $\omega_0 \in \varGamma_N\B_Y^{(d)}(\tl{U})$ is defined by
\begin{equation}
\omega_0  =\delta (y^{\prime}) dy_1 \wedge \cdots \wedge dy_d.
\end{equation}
\end{lem}

\noindent{\bf Proof of Proposition \ref{prp:7-5}}

\medskip
Let $x=(x_1,\cdots, x_n)$ be a local coordinate system of $X$ such that $M=\{x_1=\cdots =x_d=0\}$ and $(x;\xi)$ the associated coordinate system of $TX$. Set $x^{\prime}=(x_1,\cdots,x_d)$, $x^{\prime \prime}=(x_{d+1},\cdots,x_n)$ and $x=(x^{\prime},x^{\prime \prime})$.

By identifying $(or_{TM/M \utimes{X}TX}|_M)$ with $or_{M/X}$ as usual and using the isomorphism $G|_M[-n] \simeq or_X|_M$, we see that the isomorphism $\alpha$ induces isomorphisms
\begin{eqnarray}
\H_{TM}^d(\Comp_{M \utimes{X}TX} )|_M \uotimes{\Comp} (or_X|_M)
&\simeq & (or_{TM/M \utimes{X}TX}|_M) \uotimes{\Comp} (or_X|_M) \\
&\simeq & or_{M/X} \uotimes{\Comp} (or_X|_M) \\
&\simeq & or_M
\end{eqnarray}
on $M$. Via these isomorphisms, the section $u_0=1_{dx^{\prime \prime}}=1_{dx_{d+1} \wedge \cdots \wedge dx_n} \in or_M$ corresponds to the one
\begin{equation}
u_1 =[\delta (\xi^{\prime}) d\xi_1 \wedge \cdots \wedge d\xi_d ] \otimes 1_{dx_1 \wedge \cdots \wedge dx_n}
\end{equation}
of $\H_{TM}^d(\Comp_{M\utimes{X}TX})|_M \uotimes{\Comp} (or_X|_M)$ by Lemma \ref{lem:7-7}. Now let us set $L=G|_{\tLb}[-\d X] \simeq or_{\tLb}$. Then $L$ is a locally constant sheaf of rank one on $\tLb \subset M \utimes{X}TX$ whose restriction to the zero-section $M \subset \tLb$ satisfies $L|_M \simeq or_X|_M$. Moreover by the flabby resolution
\begin{equation}
0 \longrightarrow L \longrightarrow L \uotimes{\Comp} \B_{\tLb}^{(0)} \underset{d}{\longrightarrow} L \uotimes{\Comp} \B_{\tLb}^{(1)} \underset{d}{\longrightarrow} L \uotimes{\Comp} \B_{\tLb}^{(2)} \underset{d}{\longrightarrow} \cdots
\end{equation}
of $L$, we obtain an isomorphism 
\begin{eqnarray}
\lefteqn{ \H_M^d(L)|_M } \nonumber \\
&\simeq& H^d[ 0 \longrightarrow (L \uotimes{\Comp} \varGamma_M\B_{\tLb}^{(0)})|_M \underset{d}{\longrightarrow} \cdots \underset{d}{\longrightarrow} (L \uotimes{\Comp} \varGamma_M\B_{\tLb}^{(d)})|_M \longrightarrow \cdots ].
\end{eqnarray}
Therefore, a section of $\H_M^d(L)|_M$ is represented by that of
\begin{equation}
(L \uotimes{\Comp}\varGamma_M\B_{\tLb}^{(d)})|_M \simeq  (L|_M) \uotimes{\Comp}
(\varGamma_M\B_{\tLb}^{(d)}|_M) \simeq  (or_X|_M) \uotimes{\Comp} (\varGamma_M\B_{\tLb}^{(d)}|_M).
\end{equation}
Let $(\xi^{\prime}, x^{\prime \prime})$ be the coordinate system of $\tLb$ induced by that of $M \utimes{X}TX$. Then the image of the section $u_1$ by the isomorphism
\begin{equation}
\H_{TM}^d(\Comp_{M \utimes{X}TX})|_M \uotimes{\Comp} (or_X|_M) \simto \H_M^d(L)|_M
\end{equation}
induced by $\gamma$ is represented by
\begin{equation}
1_{dx_1 \wedge \cdots \wedge dx_n} \otimes \delta (\xi^{\prime}) d\xi_1 \wedge \cdots \wedge d\xi_d \in (L|_M) \uotimes{\Comp}(\varGamma_M \B_{\tLb}^{(d)}|_M).
\end{equation}
Here $1_{dx_1 \wedge \cdots \wedge dx_n}$ is considered as a section of $L|_M$ by the isomorphism $L|_M \simeq or_X|_M$.

Let us set
\begin{equation}
u_2=[1_{dx_1\wedge \cdots \wedge dx_n} \otimes \delta (\xi^{\prime}) d\xi_1 \wedge \cdots \wedge d\xi_d ] \in \H_M^d(L)|_M.
\end{equation}
Then it remains to show that the section $u_2$ is sent to
\begin{equation}
\sgn (\id-\phi^{\prime}) \cdot 1_{dx^{\prime \prime}}=\sgn(\id-\phi^{\prime}) \cdot u_0 \in or_M
\end{equation}
by $\beta$. In order to look at the isomorphism $\beta$ more precisely, consider the standard commutative diagram
\begin{equation}
\xymatrix@C=15mm{
T_{\Delta_X}(X\times X) \ar@{^{(}->}[r]^{s_0} \ar[d]^{\tau_X} & \tl{(X\times X)_{\Delta_X}} \ar[d]^{p_0} & \Omega=\{t_0>0\} \ar@{_{(}->}[l]_{j_0} \ar@{->>}[dl]^{\tl{p_0}} \\
\Delta_X \ar@{^{(}->}[r] & X \times X,}
\end{equation}
where $\tl{(X\times X)_{\Delta_X}}$ is the normal deformation of $X\times X$ along $\Delta_X$ and $t_0 \colon \tl{(X\times X)_{\Delta_X}} \longrightarrow \Real$ is the deformation parameter. Also, let $\widehat{\Omega}$ be an open subset of $\tl{(X\times X)_{\Delta_X}}$ defined by $t_0 \neq 0$ and $\widehat{p_0} \colon \widehat{\Omega} \longtwoheadrightarrow X \times X$ the restriction of $p_0$ to $\widehat{\Omega}$. Then the closure of $\widehat{p_0}^{-1}\Delta_X$ in $\tl{(X \times X)_{\Delta_X}}$ is a closed submanifold of $\tl{(X\times X)_{\Delta_X}}$, which we shall denote by $\tl{\Delta_X}$. Note that the isomorphism $(\beta|_M)^{-1}$ is the restriction of
\begin{eqnarray}
\nu_{\Delta_X}(h_*(\RG_M(\omega_X)))
&\simeq & \nu_{\Delta_X}(\RG_{\Delta_X}(h_*\omega_X)) \\
&= & {s_0}^{-1}Rj_{0*} \tl{p_0}^{-1}\RG_{\Delta_X}(h_*\omega_X) \\
&\simeq & {s_0}^{-1}\RG_{\tl{\Delta_X}}(Rj_{0*} \tl{p_0}^{-1}h_*\omega_X) \\
&\simto & \RG_{\Delta_X}(\nu_{\Delta_X}(h_*\omega_X))
\end{eqnarray}
to $M \simeq \Delta_M \subset \Delta_X$ ($\Delta_X$ is the zero-section of $T_{\Delta_X}(X\times X)$). Now let us consider also the following commutative diagram for the normal deformation $\tl{(\Gamma_{\phi})_M}$ of $\Gamma_{\phi}$ along $M \simeq \Delta_M \subset \Gamma_{\phi}$: 
\begin{equation}
\xymatrix@C=15mm{
T_M\Gamma_{\phi} \ar@{^{(}->}[r]^{s_2} \ar[d] & \tl{(\Gamma_{\phi})_M} \ar[d]_{p_2} & \Omega_{\Gamma}=\{t_2 >0\} \ar@{_{(}->}[l]_{j_2} \ar@{->>}[ld]^{\tl{p_2}} \\
M\simeq \Delta_M \ar@{^{(}->}[r] &\Gamma_{\phi},}
\end{equation}
where $t_2\colon \tl{(\Gamma_{\phi})_M} \longrightarrow \Real$ is the deformation parameter. Set $\widehat{\Omega_{\Gamma}} :=\{ t_2 \neq 0\} \subset \tl{(\Gamma_{\phi})_M}$ and $\widehat{p_2} :=p_2|_{\widehat{\Omega_{\Gamma}}}$. Then we obtain a closed submanifold $\tl{M}$ by taking the closure of $\widehat{p_2}^{-1}(M)$ in $\tl{(\Gamma_{\phi})_M}$. Since there exists an isomorphism
\begin{equation}
g^{\prime}_* \nu_M(\omega_{\Gamma_{\phi}}) \simeq \nu_{\Delta_X}(g_*\omega_{\Gamma_{\phi}}) \simeq \nu_{\Delta_X}(h_*\omega_X)
\end{equation}
(see \eqref{eq:7-3} and the Proof of Proposition \ref{prp:7-1}), the isomorphism $(\beta|_M)^{-1}$ is the restriction of
\begin{eqnarray}
\nu_M(\RG_M(\omega_{\Gamma_{\phi}}))
&\simeq & s_2^{-1}Rj_{2 *}\tl{p_2}^{-1} \RG_M(\omega_{\Gamma_{\phi}}) \label{eq:7-73}\\
&\simeq & s_2^{-1}\RG_{\tl{M}}(Rj_{2 *} \tl{p_2}^{-1}\omega_{\Gamma_{\phi}}) \\
&\simto & \RG_M (\nu_M(\omega_{\Gamma_{\phi}}))\label{eq:7-75}
\end{eqnarray}
to the zero-section $M$ of $T_M\Gamma_{\phi}$. Note that here we are identifying $T_M\Gamma_{\phi}$ with its image $\tLb \subset M \utimes{X}TX$ by $g^{\prime} \colon T_M\Gamma_{\phi} \longhookrightarrow T_{\Delta_X}(X \times X) \simeq TX$. Now let us set $L^{\prime}:=\nu_M(\omega_{\Gamma_{\phi}}) [-\d X]  \simeq or_{T_{M}\Gamma_{\phi}}$. Then $L^{\prime}$ is a locally constant sheaf of rank one on $T_M\Gamma_{\phi}$, and via the identification $T_M\Gamma_{\phi} \simeq \tLb$ we have an isomorphism $L^{\prime }\simeq L$. Note also that we have $L^{\prime}|_M \simeq L|_M \simeq or_X|_M$ by identifying $\Gamma_{\phi}$ with $X$. As a consequence, from \eqref{eq:7-73}-\eqref{eq:7-75} we obtain isomorphisms
\begin{eqnarray}
or_M
&\underset{{\bf P}}{\simto}& \H_M^d(or_{\Gamma_{\phi}})|_M \\
&\underset{{\bf Q}}{\simto}& \H_M^d(L^{\prime})|_M \\
&\underset{{\bf R}}{\simto}& \H_M^d(L)|_M
\end{eqnarray}
induced by $(\beta|_M)^{-1}$. Identify $\Gamma_{\phi}$ with $X$ as usual and use the local coordinates $x=(x^{\prime},x^{\prime \prime})$ also for $\Gamma_{\phi}$ so that we have $M \simeq \Delta_M =\{x^{\prime}=0\}$ in $\Gamma_{\phi}$. Let $(\eta^{\prime};x^{\prime \prime})$ be the associated local coordinates of $T_M\Gamma_{\phi}$. Then by the local coordinate system $(\eta^{\prime};x^{\prime \prime})$ (resp. $(\xi^{\prime};x^{\prime \prime})$ ) of $T_M\Gamma_{\phi}$ (resp. $\tLb$) the isomorphism $T_M\Gamma_{\phi} \simeq \tLb$ is given by
\begin{equation}
(\eta^{\prime};x^{\prime \prime}) \longmapsto \( \( \id-\dfrac{\partial \phi_1}{\partial x^{\prime}}\bigg|_{x^{\prime}=0}\)(\eta^{\prime});x^{\prime \prime}\),
\end{equation}
where we set $\phi(x)=(\phi_1(x),\phi_2(x))$, $\phi_1(x) \in \Real^d$, $\phi_2(x)\in \Real^{n-d}$ in the local coordinate system $x=(x^{\prime},x^{\prime \prime})$ of $X$. In particular, via the isomorphism $T_M\Gamma_{\phi} \simeq \tLb$ the differential form $\delta (\eta^{\prime}) d\eta_1 \wedge \cdots \wedge d\eta_d$ on $T_M\Gamma_{\phi}$ corresponds to the one $\sgn(\id-\phi^{\prime})\cdot \delta (\xi^{\prime}) d\xi_1 \wedge \cdots \wedge d\xi_d$ on $\tLb$. From now on, we shall determine the image of $u_0=1_{dx^{\prime \prime}} \in or_M$ by $(\beta|_M)^{-1}={\bf R} \circ {\bf Q} \circ {\bf P}$ in $\H_M^d(L)|_M$ and compare it with the previous one $u_2 \in \H_M^d(L)|_M$. First, by Lemma \ref{lem:7-6}, via the isomorphism ${\bf P}$ the section $u_0 \in or_M$ corresponds to the one of $\H_M^d(or_{\Gamma_{\phi}})|_M$ represented by
\begin{eqnarray}
1_{dx_1 \wedge \cdots \wedge dx_n} \otimes \delta (x^{\prime}) dx_1 \wedge \cdots \wedge dx_d &\in & (or_{\Gamma_{\phi}}|_M) \uotimes{\Comp}(\varGamma_M\B_{\Gamma_{\phi}}^{(d)}|_M) \\
&\simeq & (or_X|_M) \uotimes{\Comp} (\varGamma_M \B_X^{(d)}|_M).
\end{eqnarray}
It is also easy to see that via the isomorphism ${\bf Q}$ this section is sent to the one of $\H_M^d(L^{\prime})|_M$ represented by
\begin{equation}
1_{dx_1\wedge \cdots \wedge dx_n} \otimes \delta(\eta^{\prime}) d\eta_1 \wedge \cdots \wedge d\eta_d \in (L^{\prime}|_M) \uotimes{\Comp} (\varGamma_M \B_{T_M\Gamma_{\phi}}^{(d)}|_M).
\end{equation}
Here we used the isomorphism $L^{\prime}|_M \simeq or_X|_M$ to regard $1_{dx_1\wedge \cdots \wedge dx_n}$ as a section of $L^{\prime}|_M$. Finally by the isomorphism ${\bf R}$, this last section is sent to the one of $\H_M^d(L)|_M$ represented by
\begin{equation}
\sgn(\id-\phi^{\prime}) \cdot 1_{dx_1 \wedge \cdots \wedge dx_n} \otimes \delta (\xi^{\prime}) d\xi_1 \wedge \cdots \wedge d\xi_d \in (L|_M) \uotimes{M} (\varGamma_M\B_{\tLb}^{(d)}|_M),
\end{equation}
where we used $L|_M \simeq or_X|_M$. The section of $\H_M^d(L)|_M$ thus obtained differs from $u_2$ by $\sgn(\id-\phi^{\prime})$ as expected. This completes the proof. \qed

\bigskip
To conclude, we obtained the following result. 

\begin{prp}\label{prp:7-8}
In the situation as above, we have
\begin{equation}
LC(F,\Phi)=\sgn(\id-\phi^{\prime}) \cdot \overline{LC(F,\Phi)}.
\end{equation}
\end{prp}

\section{Functorial properties of Lefschetz cycles}\label{sec:8}

In this section, we study functorial properties of our Lefschetz cycles. We obtain direct and inverse image theorems for Lefschetz cycles, which extend naturally those for Kashiwara's characteristic cycles proved in \cite[Chapter IX]{K-S}.

\subsection{Direct image theorem}\label{sec:7-4}

Let $f \colon Y \longrightarrow X$ be a morphism of real analytic manifolds. Assume that we are given two morphisms $\phi_X \colon X \longrightarrow X$ and $\phi_Y \colon Y \longrightarrow Y$ such that the diagram
\begin{equation}
\xymatrix{
Y \ar[r]^f \ar[d]_{\phi_Y} & X \ar[d]^{\phi_X} \\
Y \ar[r]^f & X}
\end{equation}
commutes. Consider also an object $G$ of $\Dc(Y)$ such that $f$ is proper on $\supp(G)$ and a morphism
\begin{equation}
\Phi_Y \colon \phi_Y^{-1}G \longrightarrow G
\end{equation}
in $\Dc(Y)$. Then $Rf_*G \in \Dc(X)$ and we obtain a natural morphism
\begin{equation}
\Phi_X \colon \phi_X^{-1}Rf_*G \longrightarrow Rf_*G
\end{equation}
induced by $\Phi_Y$. Our aim in this subsection is to compare the Lefschetz cycle of $(G,\Phi_Y)$ with that of $(Rf_*G, \Phi_X)$. Let $M$ be a smooth fixed point component of $\phi_X$ such that $f(\supp(G)) \cap M$ is compact. Also let $\{N_i\}_{i \in I}$ be the set of all fixed point components $N_i$ of $\phi_Y$ such that $N_i \subset f^{-1}(M)$ and $\supp (G) \cap N_i \not= \emptyset$. Note that $I$ is a finite set by our assumptions. Set $N:=\bigsqcup_{i\in I}N_i$ and assume that $N$ is smooth. We also assume that $\Gamma_{\phi_X} \subset X\times X$ (resp. $\Gamma_{\phi_Y} \subset Y \times Y$) intersects with $\Delta_X$ in $X \times X$ (resp. $\Delta_Y$ in $Y \times Y$) cleanly along $M$ (resp. $N$) as in previous sections. For the sake of simplicity, denote $M\utimes{X} \{ T_{\Gamma_{\phi_X}}^*(X\times X) \cap T^*_{\Delta_X}(X \times X) \} \simeq T^*M$, $N \utimes{Y} \{ T_{\Gamma_{\phi_Y}}^*(Y\times Y)\cap T^*_{\Delta_Y}(Y \times Y) \} \simeq T^*N$ simply by $\Lb$, $\sLb$ respectively. Then we obtain two Lefschetz bundles
\begin{gather}
\Lb \subset T_{\Gamma_{\phi_X}}^*(X\times X) \cap T^*_{\Delta_X}(X \times X), \\
\sLb \subset T_{\Gamma_{\phi_Y}}^*(Y\times Y) \cap T^*_{\Delta_Y}(Y \times Y)
\end{gather}
and the Lefschetz cycles
\begin{gather}
LC(G,\Phi_Y)_N \in H_{\SS(G) \cap \sLb}^0(\sLb; \pi_N^{-1}\omega_N), \\
LC(Rf_*G, \Phi_X)_M \in H_{\SS(Rf_*G) \cap \Lb}^0(\Lb ; \pi_M^{-1}\omega_M)
\end{gather}
in them. Note that by setting $\sLb_i:= N_i \times_{N} \sLb$ we have the direct sum decompositions $\sLb =\bigsqcup_{i\in I}\sLb_i \simeq \bigsqcup_{i\in I}T^*N_i$ and 
\begin{equation} 
LC(G,\Phi_Y)_N=\dsum_{i\in I} LC(G,\Phi_Y)_{N_i}, 
\end{equation}  
where $LC(G,\Phi_Y)_{N_i} \in H_{\SS(G) \cap \sLb_i}^0(\sLb_i; \pi_{N_i}^{-1}\omega_{N_i})$. 
Now let us set $g =f|_N \colon N =\bigsqcup_{i\in I}N_i \longrightarrow M$ and consider the natural morphisms
\begin{equation}
T^*N \overset{\t{g}}{\longleftarrow} N \utimes{M}T^*M \overset{g_{\pi}}{\longrightarrow} T^*M
\end{equation}
induced by $g$. Take a closed conic subanalytic Lagrangian subset $\Lambda =\bigsqcup_{i \in I}\Lambda_i$ of $T^*N =\bigsqcup_{i \in I}T^*N_i$ such that $\SS(G) \cap \sLb \subset \Lambda$ and $g_{\pi}$ is proper on $\t{g}^{-1}(\Lambda)$. Set $\Lambda^{\prime} =\t{g}^{-1}(\Lambda)$ and $\Lambda^{\prime \prime}=g_{\pi}(\Lambda^{\prime})$. Then there exists a morphism
\begin{equation}\label{eq:7-94}
g_* \colon H_{\Lambda}^0(T^*N ;\pi_N^{-1}\omega_N) \longrightarrow H_{\Lambda^{\prime \prime}}^0(T^*M;\pi_M^{-1}\omega_M)
\end{equation}
of Lagrangian cycles induced by $g$ (see \cite[Proposition 9.3.2 (i)]{K-S}).

Note that by the commutativity of the diagram 
\begin{equation}
\xymatrix{
Y \ar[r]^f& X \\
N \ar@{^{(}->}[u]_{i_N} \ar[r]^{g} & M \ar@{^{(}->}[u]_{i_M}}
\end{equation}
we have a commutative diagram 
\begin{equation}
\xymatrix@R=5mm@C=15mm{
N \utimes{M} \Lb \ar[r]^{ \rho_0 } \ar@{-}[d]^{\wr} & \sLb \ar@{-}[d]^{\wr} \\
N \utimes{M}T^*M  \ar[r]^{ \t{g} } & T^*N,}
\end{equation}
where $\rho_0$ is the restriction of the natural morphism $\t{f} \colon Y \utimes{X} T^*X \longrightarrow T^*Y$ to $N \utimes{M} \Lb \subset Y \utimes{X} T^*X$.

\begin{thm}\label{thm:7-8}
In the situation as above, we have
\begin{equation}
LC(Rf_*G, \Phi_X)_M =g_*(LC(G, \Phi_Y)_N)
\end{equation}
in $T^*M$. More precisely, for the morphism
\begin{equation}
(g_i)_* \colon H_{\Lambda_i}^0(T^*N_i ;\pi_{N_i}^{-1} \omega_{N_i}) \longrightarrow H_{\Lambda^{\prime \prime}}^0 (T^*M;\pi_M^{-1}\omega_M)
\end{equation}
of Lagrangian cycles induced by $g_i =f|_{N_i} \colon N_i \longrightarrow M$ we have 
\begin{equation}
LC(Rf_*G, \Phi_X)_M = \dsum_{i \in I} (g_i)_*(LC(G, \Phi_Y)_{N_i}). 
\end{equation}
\end{thm}

\begin{proof}
The proof is similar to that of \cite[Proposition 9.4.2]{K-S}. Let $\delta \colon Y \longhookrightarrow X\times Y$ be a morphism defined by $y \longmapsto (f(y),y)$. Then the image of $\delta$ is the graph $\Delta:=\Gamma_f \subset X\times Y$ of $f$. Let $\delta_X \colon X \longhookrightarrow X \times X$ and $\delta_Y \colon Y \longhookrightarrow Y \times Y$ be the diagonal embeddings of $X$ and $Y$ respectively. Then we obtain a commutative diagram
\begin{equation}
\xymatrix@R=3mm@C=3mm{
Y \times Y \ar[rr]^{f_1} & & X \times Y \ar[rr]^{f_2} & & X \times X \\
 & & & \Box & \\
Y \ar@{_{(}->}[uu]^{\delta_Y} \ar@{=}[rr]^{\id_Y} & & Y \ar@{_{(}->}[uu]^{\delta} \ar[rr]^f & & X, \ar@{_{(}->}[uu]^{\delta_X}}
\end{equation}
where we set $f_1 :=f \times \id_Y$ and $f_2 :=\id_X \times f$. We also need the commutative diagram
\begin{equation}
\xymatrix@R=3mm@C=3mm{
Y \times Y \ar[rr]^{f_1} & & X \times Y \ar[rr]^{f_2} & & X \times X \\
 & & & \Box & \\
Y \ar@{_{(}->}[uu]^{h_Y} \ar@{=}[rr]^{\id_Y} & & Y \ar@{_{(}->}[uu]^{h} \ar[rr]^f & & X, \ar@{_{(}->}[uu]^{h_X}}
\end{equation}
where the morphisms $h_X$, $h_Y$ and $h$ are defined by $x \longmapsto (\phi_X(x),x)$, $y \longmapsto (\phi_Y(y),y)$ and $y \longmapsto ((\phi_X \circ f)(y),y)=((f \circ \phi_Y)(y),y)$ respectively. Set $\Gamma:=h(Y) \subset X \times Y$. Since the morphism $f$ is decomposed as 
\begin{equation}
\xymatrix@C=15mm{
Y \ar@{^{(}->}[r]^{\delta} \ar[dr]_f & X\times Y \ar[d]^{p_X} \\ & X}
\end{equation}
and there exists a commutative diagram
\begin{equation}
\xymatrix@C=15mm{
Y \ar@{^{(}->}[r]^{\delta} \ar[d]^{\phi_Y} & X \times Y \ar[r]^{p_X} \ar[d]^{\phi_X \times \phi_Y} & X \ar[d]^{\phi_X} \\
Y \ar@{^{(}->}[r]^{\delta} & X\times Y \ar[r]^{p_X} & X,}
\end{equation}
we may reduce the problem to the case of closed embeddings and that of smooth maps. 

First, let us consider the case where $f \colon Y \longrightarrow X$ is smooth and proper on $\supp(G)$. Then we have the following lemma whose proof is similar to that of Lemma \ref{lem:4-5} (use also the proof of Proposition \ref{prp:7-1}). 

\begin{lem}\label{lem:7-9}
Assume that $f$ is smooth. Then 
\begin{enumerate}
\item By identifying $\Delta$ with $Y$, we have
\begin{equation}
\Gamma \cap \Delta \simeq f^{-1}(M).
\end{equation}
\item By identifying $T^*_{\Delta}(X \times Y)$ with $Y \utimes{X}T^*X$, we have
\begin{equation}
T_{\Gamma}^*(X \times Y) \cap T_{\Delta}^*(X\times Y) \simeq f^{-1}(M)
\utimes{M}\Lb.
\end{equation}
\item The support of $\mu_{\Delta}(h_*\omega_Y)$ is $f^{-1}(M) \utimes{M}\Lb \subset Y \utimes{X}T^*X$ and there exists an isomorphism 
\begin{equation}
\mu_{\Delta}(h_*\omega_Y)|_{f^{-1}(M) \utimes{M}\Lb} \simeq \pi_0^{-1}\omega_{f^{-1}(M)}, 
\end{equation}
where $\pi_0 \colon f^{-1}(M) \utimes{M} \Lb \simeq f^{-1}(M) \utimes{M}T^*M \longrightarrow f^{-1}(M)$ is the projection. 
\end{enumerate}
\end{lem}
Now consider also the natural morphisms
\begin{equation}
T^*Y \overset{\t{f}}{\longleftarrow} Y \utimes{X}T^*X \overset{f_{\pi}}{\longrightarrow} T^*X
\end{equation}
induced by $f$. For short, let us set $S=\SS(G)$, $S^{\prime}=\t{f}^{-1}(S)$ and $S^{\prime \prime}=f_{\pi}(S^{\prime})$. Then by Lemma \ref{lem:7-9} and the morphisms of functors
\begin{gather}
\t{f_1}^{-1} \circ \mu_{\Delta_Y} \longrightarrow \mu_{\Delta}\circ Rf_{1!}, \\
Rf_{2 \pi !} \circ \mu_{\Delta} \longrightarrow \mu_{\Delta_X} \circ Rf_{2!}
\end{gather}
obtained by \cite[Proposition 4.3.4]{K-S} we obtain the commutative diagram Diagram 8.a (we omit the symbols $R$ of right derived functors etc. to simplify the notation). Here the middle vertical arrows are induced by
\begin{equation}
h_*(f^{-1}f_*G \otimes \D G) \longrightarrow  h_*(G \otimes \D G) \longrightarrow  h_*\omega_Y.
\end{equation}

\newpage
\begin{center}{\small
\rotatebox[origin=c]{90}{
$
\xymatrix{
{\rm Hom}(G,G) \ar[rr]& & {\rm Hom}(f_*G, f_*G) \\
\RG_{\Delta_Y}(Y \times Y; G \boxtimes \D G) \ar[u]^{\wr} \ar@{-}[d]^{\wr} \ar[r]& \RG_{\Delta}(X\times Y; f_*G \boxtimes \D G) \ar[r] \ar@{-}[d]^{\wr} & \RG_{\Delta_X}(X \times X; f_*G \boxtimes \D f_*G) \ar@{-}[d]^{\wr} \ar[u]^{\wr} \\
\RG_S(T^*Y; \mu_{\Delta_Y}(G \boxtimes \D G)) \ar[r] \ar[d] & \RG_{S^{\prime}}(Y \utimes{X}T^*X; \mu_{\Delta}(f_*G \boxtimes \D G)) \ar[r] \ar[d] & \RG_{S^{\prime \prime}}(T^*X; \mu_{\Delta_X}(f_*G \boxtimes \D f_*G)) \ar[d] \\
%\RG_S(T_*Y;\mu_{\Delta_Y}(h_{Y*}h_Y^{-1}(G \boxtimes \D G))) \ar[r] \ar[d] & \RG_{S^{\prime}}(Y \utimes{X}T^*X;\mu_{\Delta}(h_*h^{-1}(f_*G \boxtimes \D G))) \ar[r] \ar[d] & \RG_{S^{\prime \prime}}(T^*X; \mu_{\Delta_X}(h_{X*}h_X^{-1}(f_*G \boxtimes \D f_*G))) \ar[d] \\
\RG_S(T_*Y;\mu_{\Delta_Y}(h_{Y*}(\phi_Y^{-1}G \otimes \D G))) \ar[d]^{\Phi_Y} & \RG_{S^{\prime}}(Y \utimes{X}T^*X;\mu_{\Delta}(h_*(f^{-1}\phi_X^{-1}f_*G \otimes \D G))) \ar[r] \ar[d]^{\Phi_X} & \RG_{S^{\prime \prime}}(T^*X; \mu_{\Delta_X}(h_{X*}(\phi_X^{-1}f_*G \otimes \D f_*G))) \ar[d]^{\Phi_X} \\
\RG_S(T_*Y;\mu_{\Delta_Y}(h_{Y*}(G \otimes \D G))) \ar[d] & \RG_{S^{\prime}}(Y \utimes{X}T^*X;\mu_{\Delta}(h_*(f^{-1}f_*G \otimes \D G))) \ar[r] \ar[d]& \RG_{S^{\prime \prime}}(T^*X; \mu_{\Delta_X}(h_{X*}(f_*G \otimes \D f_*G))) \ar[d] \\
\RG_S(T_*Y;\mu_{\Delta_Y}(h_{Y*}\omega_Y))\ar[r] \ar[d] & \RG_{S^{\prime}}(Y \utimes{X}T^*X;\mu_{\Delta}(h_*\omega_Y)) \ar[r] \ar[d] & \RG_{S^{\prime \prime}}(T^*X; \mu_{\Delta_X}(h_{X*}\omega_X)) \ar[d] \\
\RG_{S\cap \sLb}(T_*N; \pi_N^{-1}\omega_N) \ar[r]^{\hspace*{-20mm}{\bf A}} & \RG_{S^{\prime}\cap (f^{-1}(M) \utimes{M}\Lb)}(f^{-1}(M) \utimes{M}T^*M;\pi_0^{-1}\omega_{f^{-1}(M)}) \ar[r]^{\hspace*{15mm}{\bf B}} & \RG_{S^{\prime \prime} \cap \Lb}(T^*M; \pi_M^{-1}\omega_M).}
$}}

\bigskip
Diagram 8.a
\end{center}

\newpage
\noindent
Consider the commutative diagram
\begin{equation}
\xymatrix@R=3mm@C=6mm{
\sLb \simeq T^*N \ar[ddrr]_{\pi_N} & &N \utimes{M}\Lb \ar@{_{(}->}[ll]_{\t{g}} \ar[dd]^{\pi} \ar@{^{(}->}[rr]^{k} & & f^{-1}(M) \utimes{M}\Lb \ar[dd]^{\pi_0}\\
 & & & \Box & \\
 & &N \ar@{^{(}->}[rr] & & f^{-1}(M),}
\end{equation}
where $\pi \colon N \utimes{M}\Lb \longrightarrow N$ (resp. $k \colon N \utimes{M}\Lb \longhookrightarrow f^{-1}(M) \utimes{M}\Lb$) is the projection (resp. the inclusion map). Then the morphism ${\bf A}$ is decomposed as follows.
\begin{eqnarray}
\RG_{S \cap \sLb}(\sLb;\pi_N^{-1}\omega_N)
&\overset{\t{g}^{-1}}{\longrightarrow} & \RG_{\t{g}^{-1}(S \cap \sLb)} (N \utimes{M}\Lb ;\pi^{-1}\omega_N) \\
&\overset{\int_k}{\longrightarrow} & \RG_{S^{\prime} \cap (f^{-1}(M) \utimes{M}\Lb)} (f^{-1}(M) \utimes{M} \Lb ; \pi_0^{-1}\omega_{f^{-1}(M)}).
\end{eqnarray}
Moreover we can easily show that the morphism ${\bf B}$ is induced by the topological integration morphism
\begin{equation}
\dint_q \colon Rq_!q^! \omega_M \simeq Rq_! \omega_{f^{-1}(M)} \longrightarrow \omega_M
\end{equation}
for $q \colon f^{-1}(M) \longrightarrow M$. Therefore the composite of ${\bf A}$ and ${\bf B}$ is factorized as follows
\begin{equation}\label{eq:8-25} 
\xymatrix@C=20mm{
\RG_{S \cap \sLb}(\sLb; \pi_N^{-1}\omega_N) \ar[r]^{\t{g}^{-1}} \ar[rd]_{{\bf B} \circ {\bf A}} & \RG_{\t{g}^{-1}(S \cap \sLb)}(N \utimes{M}\Lb;\pi^{-1}\omega_N) \ar[d]^{\int_g} \\
& \RG_{S^{\prime \prime} \cap \Lb}(\Lb; \pi_M^{-1}\omega_M).}
\end{equation}
This implies that the morphism ${\bf B}\circ {\bf A}$ is the push-forward of Lagrangian cycles \eqref{eq:7-94}. 

Next, consider the case where $f \colon Y \longrightarrow X$ is a closed embedding. Then similarly we have the following. 

\begin{lem}\label{lem:8-3}
Assume that $f$ is a closed embedding. Then 
\begin{enumerate}
\item By identifying $\Delta$ with $Y$, we have
\begin{equation}
\Gamma \cap \Delta \simeq f^{-1}(M)=N.
\end{equation}
\item By identifying $T^*_{\Delta}(X \times Y)$ with $Y \utimes{X}T^*X$, we have
\begin{equation}
T_{\Gamma}^*(X \times Y) \cap T_{\Delta}^*(X\times Y) \simeq \t{f}^{-1} \sLb
\end{equation}
and $\t{f}^{-1} \sLb \supset f_{\pi}^{-1} \Lb \simeq N \utimes{M} \Lb$.
\item The support of $\mu_{\Delta}(h_*\omega_Y)$ is $\t{f}^{-1} \sLb \subset Y \utimes{X}T^*X$ and there exists an isomorphism 
\begin{equation}
\mu_{\Delta}(h_*\omega_Y)|_{\t{f}^{-1} \sLb} \simeq \pi_1^{-1}\omega_{N}, 
\end{equation}
where $\pi_1 \colon \t{f}^{-1} \sLb \simeq \t{f}^{-1} T^*N \longrightarrow N$ is the projection. 
\end{enumerate}
\end{lem}

Now, by Lemma \ref{lem:8-3} we obtain a new commutative diagram by replacing the bottom horizontal arrows in Diagram 8.a by 
\begin{eqnarray}
\RG_{S \cap \sLb}(T^*N ;\pi_N^{-1}\omega_N) & \overset{{\bf C}}{\longrightarrow} & \RG_{S^{\prime}\cap \t{f}^{-1} \sLb} ( \t{f}^{-1} \sLb  ;\pi_1^{-1}\omega_{N}) \\
& \overset{{\bf D}}{\longrightarrow} & \RG_{S^{\prime \prime} \cap \Lb}(T^*M; \pi_M^{-1}\omega_M). 
\end{eqnarray}
Here the morphism ${\bf C}$ is induced by $\id \longrightarrow R \t{f}_* \circ \t{f}^{-1}$. Moreover the morphism ${\bf D}$ is decomposed as
\begin{eqnarray}
\RG_{S^{\prime}\cap \t{f}^{-1} \sLb} ( \t{f}^{-1} \sLb  ;\pi_1^{-1}\omega_{N})
& \overset{\alpha}{\longrightarrow} & \RG_{S^{\prime}\cap (N \utimes{M} \Lb)} (N \utimes{M}T^*M ;\pi^{-1}\omega_{N}) \\
& \overset{\beta}{\longrightarrow} & \RG_{S^{\prime \prime} \cap \Lb}(T^*M  ; \pi_M^{-1}\omega_M),
\end{eqnarray}
where the morphism $\alpha$ (resp. $\beta$) is induced by the restriction to $N \utimes{M} T^*M \subset \t{f}^{-1} \sLb$ (resp. by the natural morphism $\int_g  \colon g_* \omega_N \simeq \RG_{N}\omega_M \longrightarrow \omega_M$). Hence the composite of ${\bf C}$ and ${\bf D}$ is factorized as in the diagram \eqref{eq:8-25}.
This completes the proof.
\qed
\end{proof}

Since via the characteristic cycle maps $CC$ (see Proposition \ref{prp:2-10}) the push-forward $g_*$ of Lagrangian cycles corresponds to the topological integral 
\begin{equation} 
\int_g \colon \CF(N)_{\Comp} \longrightarrow \CF(M)_{\Comp}
\end{equation} 
of constructible functions, by Theorem \ref{thm:5-1} we obtain the following result.

\begin{cor}(see also {\rm \cite{M-T-3}}) For the local contributions $c(G, \Phi_Y)_{N_i}$ and $c(Rf_*G, \Phi_X)_M$, we have 
\begin{equation}
c(Rf_*G, \Phi_X)_M = \dsum_{i \in I} c(G, \Phi_Y)_{N_i}.
\end{equation}
\end{cor}

\subsection{Inverse image theorem}

In this subsection, we establish the inverse image theorem for Lefschetz cycles. We mainly inherit the notations and the situation treated in Section \ref{sec:7-4}. However, here $M$ and $N$ are smooth fixed point components of $\phi_X$ and $\phi_Y$ respectively satisfying just the condition $f(N) \subset M$. Consider an object $F$ of $\Dc(X)$ and a morphism
\begin{equation}
\Phi_X \colon \phi_X^{-1}F \longrightarrow F
\end{equation}
in $\Dc(X)$. Then $f^{-1}F \in \Dc(Y)$ and we obtain a natural morphism
\begin{equation}
\Phi_Y \colon \phi_Y^{-1}f^{-1}F \longrightarrow f^{-1}F
\end{equation}
induced by $\Phi_X$. Assuming the same conditions on $\phi_X$, $\phi_Y$ etc. and keeping the same notations for $\Lb$, $\sLb$ etc. as in Section \ref{sec:7-4}, we obtain the Lefschetz cycles
\begin{gather}
LC(F,\Phi_X)_M \in H_{\SS(F) \cap \Lb}^0(\Lb; \pi_M^{-1}\omega_M), \\
LC(f^{-1}F, \Phi_Y)_N \in H_{\SS(f^{-1}F) \cap \sLb}^0(\sLb ;\pi_N^{-1}\omega_N).
\end{gather}
Set $g =f|_N \colon N \longrightarrow M$ as before and consider the natural morphisms
\begin{equation}
T^*N \overset{\t{g}}{\longleftarrow} N \utimes{M}T^*M \overset{g_{\pi}}{\longrightarrow} T^*M
\end{equation}
induced by $g$. Let $\Lambda \subset \SS(F) \cap \Lb \subset \Lb \simeq T^*M$ be the support of $LC(F,\Phi_X)_M$ and set $\Lambda^{\prime}=g_{\pi}^{-1}(\Lambda)$ and $\Lambda^{\prime\prime}=\t{g}(\Lambda^{\prime})$. If $\t{g}$ is proper on $\Lambda^{\prime}$ (e.g. if $f$ is non-characteristic for $F$ on an open neighborhood of $N$), then there exists a morphism
\begin{equation}\label{eq:8-42}
g^* \colon H_{\Lambda}^0(T^*M ;\pi_M^{-1}\omega_M) \longrightarrow H_{\Lambda^{\prime \prime}}^0(T^*N;\pi_N^{-1}\omega_N)
\end{equation}
of Lagrangian cycles induced by $g$ (see \cite[Proposition 9.3.2 (ii)]{K-S}).

\begin{thm}\label{thm:7-10}
In the situation as above, assume moreover that $f$ is non-characteristic for $F$ on an open neighborhood of $N$. Then we have
\begin{equation}
LC(f^{-1}F, \Phi_Y)_N =\sgn(\id-\phi_X^{\prime})\cdot \sgn(\id-\phi_Y^{\prime}) \cdot g^*(LC(F, \Phi_X)_M)
\end{equation}
in $T^*N$, where $\sgn(\id-\phi_X^{\prime})=\pm 1$ (resp. $\sgn(\id-\phi_Y^{\prime})=\pm 1$) is the sign of the determinant of $\id-\phi_X^{\prime} \colon T_MX \longrightarrow T_MX$ (resp. $\id-\phi_Y^{\prime} \colon T_NY \longrightarrow T_NY$).
\end{thm}

\begin{proof}
The proof is similar to that of \cite[Proposition 9.4.3]{K-S}. By Proposition \ref{prp:7-8}, it suffices to show that
\begin{equation}\label{eq:7-113}
\overline{LC(f^{-1}F,\Phi_Y)_N} =g^* (\overline{LC(F, \Phi_X)_M}).
\end{equation}
We use almost the same notations as in the proof of Theorem \ref{thm:7-8} except the ones for $\delta$ and $h$ etc. In particular, here we define a morphism $\delta \colon Y \longhookrightarrow Y\times X$ by $y \longmapsto (y,f(y))$ and set $\Delta:=\delta(Y)\subset Y \times X$. Also we define a morphism $h \colon Y \longrightarrow Y \times X$ by $y \longmapsto (\phi_Y(y),f(y))$ and set $\Gamma:=h(Y) \subset Y \times X$. Then we have the following commutative diagrams.
\begin{gather}
\xymatrix@R=3mm@C=3mm{
Y \times Y \ar[rr]^{f_2} & & Y \times X \ar[rr]^{f_1} & & X \times X \\
 & & & \Box & \\
Y \ar@{_{(}->}[uu]^{\delta_Y} \ar@{=}[rr]^{\id_Y} & & Y \ar@{_{(}->}[uu]^{\delta} \ar[rr]^f & & X, \ar@{_{(}->}[uu]^{\delta_X}}\\
\label{diag:8-53}\xymatrix{
Y \times Y \ar[r]^{f_2} & Y \times X \ar[r]^{f_1} & X \times X \\
Y \ar@{_{(}->}[u]^{h_Y} \ar@{=}[r]^{\id_Y} & Y \ar[u]^{h} \ar[r]^f & X. \ar@{_{(}->}[u]^{h_X}}
\end{gather}
Note that $h$ is not always injective. Now by using the natural morphisms
\begin{equation}
T^*Y \overset{\t{f}}{\longleftarrow} Y \utimes{X}T^*X \overset{f_{\pi}}{\longrightarrow} T^*X
\end{equation}
induced by $f$, set for short $S:=\SS(F)$, $S^{\prime}:=f_{\pi}^{-1}(S)$ and $S^{\prime \prime}:=\t{f}(S^{\prime})$. Then by the morphisms of functors
\begin{gather}
f_{1\pi}^{-1} \circ \mu_{\Delta_X} \longrightarrow \mu_{\Delta}\circ f_1^{-1}, \\
R \t{f_2}_! \circ \mu_{\Delta} \longrightarrow \mu_{\Delta_Y}(\omega_{Y/X} \otimes f_2^{-1}),
\end{gather}
obtained by \cite[Proposition 4.3.5]{K-S} we obtain the commutative diagram (similar to \cite[Diagram 9.4.6]{K-S}) Diagram 8.b. Here we omit the symbols $R$ of right derived functors etc. to simplify the notation.

\newpage
\begin{center}{\small
\rotatebox[origin=c]{90}{
$
\xymatrix{
{\rm Hom}(F,F) \ar[rr] & & {\rm Hom}(f^{-1}F,f^{-1}F) \\
\RG_{\Delta_X}(X\times X;F \boxtimes \D F) \ar[r] \ar[u]^{\wr} \ar@{-}[d]^{\wr} & \RG_{\Delta}(Y \times X; f^{-1}F \boxtimes \D F) \ar@{-}[r]^{\sim} \ar@{-}[d]^{\wr} & \RG_{\Delta_Y}(Y \times Y;f^{-1}F \boxtimes \D f^{-1}F) \ar[u]^{\wr} \ar@{-}[d]^{\wr} \\
\RG_S(T^*X; \mu_{\Delta_X}(F \boxtimes \D F)) \ar[r] \ar[d]& \RG_{S^{\prime}}(Y \utimes{X}T^*X; \mu_{\Delta}(f^{-1}F \boxtimes \D F)) \ar[r] \ar[d]& \RG_{S^{\prime \prime}}(T^*Y:\mu_{\Delta_Y}(f^{-1} F\boxtimes \D f^{-1}F)) \ar[d]\\
\RG_S(T^*X; \mu_{\Delta_X}(h_{X*}(\phi_X^{-1}F \otimes \D F))) \ar[r] \ar[d]^{\Phi_X}& \RG_{S^{\prime}}(Y \utimes{X}T^*X; \mu_{\Delta}(h_*(\phi_Y^{-1}f^{-1}F \otimes f^{-1}\D F))) \ar[r] \ar[d]^{\Phi_Y}& \RG_{S^{\prime \prime}}(T^*Y;\mu_{\Delta_Y}(h_{Y*}(\phi_Y^{-1}f^{-1} F\otimes \D f^{-1}F))) \ar[d]^{\Phi_Y}\\
\RG_S(T^*X; \mu_{\Delta_X}(h_{X*}(F \otimes \D F))) \ar[r] \ar[d]& \RG_{S^{\prime}}(Y \utimes{X}T^*X; \mu_{\Delta}(h_*(f^{-1}F \otimes f^{-1}\D F))) \ar[r] \ar[d]& \RG_{S^{\prime \prime}}(T^*Y;\mu_{\Delta_Y}(h_{Y*}(f^{-1} F\otimes \D f^{-1}F))) \ar[d]\\
\RG_S(T^*X; \mu_{\Delta_X}(h_{X*}\omega_X)) \ar[r]^{{\bf A}} \ar[d]^{\wr} & \RG_{S^{\prime}}(Y \utimes{X}T^*X; \mu_{\Delta}(h_*f^{-1}\omega_X)) \ar[r]^{{\bf B}} & \RG_{S^{\prime \prime}}(T^*Y;\mu_{\Delta_Y}(h_{Y*}\omega_Y)) \ar[d]^{\wr} \\
\RG_{S \cap \Lb}(T^*M;\pi_M^{-1}\omega_M) \ar[rr] & & \RG_{S^{\prime \prime} \cap \sLb}(T^*N;\pi_N^{-1}\omega_N).}
$}}

\bigskip
Diagram 8.b
\end{center}

\newpage
\noindent
Let us explain the construction of the bottom square in Diagram 8.b. We consider the commutative diagram:
\begin{equation}\label{eq:8-58}
\xymatrix{
Y \ar[r]^f& X \\
N \ar@{^{(}->}[u]_{i_N} \ar[r]^{g} & M \ar@{^{(}->}[u]_{i_M}.}
\end{equation}
Let
\begin{equation}
\Theta_f \colon \RG_S(T^*X;\mu_{\Delta_X}(h_{X*}\omega_X)) \longrightarrow \RG_{S^{\prime \prime}}(T^*Y;\mu_{\Delta_Y}(h_{Y*}\omega_Y))
\end{equation}
be a morphism obtained by
\begin{eqnarray}
\RG_S(T^*X;\mu_{\Delta_X}(h_{X*}\omega_X))
&\longrightarrow & \RG_{S^{\prime}}(Y \utimes{X}T^*X;f_{1\pi}^{-1}\mu_{\Delta_X}(h_{X*}\omega_X)) \\
&\longrightarrow & \RG_{S^{\prime}}(Y \utimes{X}T^*X;\mu_{\Delta}(f_1^{-1} h_{X*} \omega_X)) \\
&\longrightarrow & \RG_{S^{\prime}}(Y \utimes{X}T^*X;\mu_{\Delta}(Rh_*f^{-1}\omega_X)) \\
&\longrightarrow & \RG_{S^{\prime \prime}}(T^*Y;R \t{f_2}_! \mu_{\Delta}(Rh_*f^{-1}\omega_X)) \\
&\longrightarrow & \RG_{S^{\prime \prime}}(T^*Y;\mu_{\Delta_Y}(\omega_{Y/X} \otimes f_2^{-1}Rh_*f^{-1}\omega_X)) \\
&\longrightarrow & \RG_{S^{\prime \prime}}(T^*Y;\mu_{\Delta_Y}(h_{Y*}\omega_Y)),
\end{eqnarray}
where we used \cite[Proposition 4.3.5]{K-S} to construct the second and fifth morphisms. Moreover we used the assumption that $f$ is non-characteristic for $F$ to construct the forth one. The morphism $\Theta_f$ is the composite of ${\bf A}$ and ${\bf B}$ in Diagram 8.b. Let $\delta_g \colon N \longrightarrow N \times M$ be a morphism defined by $y \longmapsto (y, g(y))$ and note that $h_X|_M=\delta_M$, $h_Y|_N=\delta_N$ and $h|_N=\delta_g$. If we paraphrase the construction of $\Theta_f$ by using the commutative diagram
\begin{equation}
\xymatrix{
N \times N\ar[r]^{g_2} & N \times M \ar[r]^{g_1} & M \times M \\
N \ar@{_{(}->}[u]^{\delta_N} \ar@{=}[r]^{\id_N} & N \ar[u]^{\delta_g} \ar[r]^g & M, \ar@{_{(}->}[u]^{\delta_M}}
\end{equation}
instead of the one \eqref{diag:8-53}, we obtain a morphism 
\begin{equation}
\Theta_g \colon \RG_{S \cap \Lb}(T^*M;\mu_{\Delta_M}(\delta_{M*}\omega_M)) \longrightarrow \RG_{\t{g}g_{\pi}^{-1}(S \cap \Lb)}(T^*N;\mu_{\Delta_N}(\delta_{N*}\omega_N)).
\end{equation}
By the proof of \cite[Proposition 9.3.2 (ii) and Proposition 9.4.3]{K-S}, there exists a commutative diagram 
\begin{equation}\label{diag:8-43}
\xymatrix{
\RG_{S \cap \Lb}(T^*M;\mu_{\Delta_M} (\delta_{M*}\omega_M)) \ar[r]^{\Theta_g} \ar@{-}[d]^{\wr} & \RG_{\t{g}g_{\pi}^{-1}(S \cap \Lb)}(T^*N;\mu_{\Delta_N}(\delta_{N*}\omega_N)) \ar@{-}[d]^{\wr} \\
\RG_{S \cap \Lb}(T^*M;\pi_M^{-1}\omega_M) \ar[r]^{g^*} & \RG_{\t{g}g_{\pi}^{-1}(S \cap \Lb)}(T^*N;\pi_N^{-1}\omega_N). }
\end{equation}
Since we have $\t{g}g_{\pi}^{-1}(S \cap \Lb) \subset S^{\prime \prime} \cap \sLb$, by $\Theta_g$ we obtain also a morphism
\begin{equation}
\RG_{S \cap \Lb} (T^*M;\mu_{\Delta_M}(\delta_{M*}\omega_M)) \longrightarrow \RG_{S^{\prime \prime} \cap \sLb}(T^*N;\mu_{\Delta_N}(\delta_{N*}\omega_N)). 
\end{equation}
We still denote it by $\Theta_g$. In the same way, we can construct also morphisms
\begin{gather}
\Theta_{i_M} \colon \RG_{S}(T^*X;\mu_{\Delta_X}(h_{X*}\omega_X)) \longrightarrow \RG_{S\cap \Lb}(T^*M; \mu_{\Delta_M}(\delta_{M*}\omega_M)), \\
\Theta_{i_N} \colon \RG_{S^{\prime\prime}}(T^*Y;\mu_{\Delta_Y}(h_{Y*}\omega_Y)) \longrightarrow \RG_{S^{\prime\prime} \cap \sLb}(T^*N;\mu_{\Delta_N}(\delta_{N*}\omega_N)),
\end{gather}
where we used the fact that the support of $\mu_{\Delta_X}(h_{X*}\omega_X)$ (resp. $\mu_{\Delta_Y}(h_{Y*}\omega_Y)$) is contained in $\Lb$ (resp. $\sLb$) and $\t{i_M}$ (resp. $\t{i_N}$) is proper on $S \cap \Lb$ (resp. $S^{\prime \prime} \cap \sLb$). Now recall that the morphism $f$ is the composite of the graph embedding $\delta \colon Y \longhookrightarrow Y \times X$ and the projection $p_X \colon Y \times X \longtwoheadrightarrow X$. Then we may assume that $f$ is a graph embedding or a projection. In both cases, since the constructions of the morphisms $\Theta_f$, $\Theta_g$, $\Theta_{i_M}$ and $\Theta_{i_N}$ are similar, we obtain the following commutative diagram:
\begin{equation}
\xymatrix@R=8mm@C=25mm{
\RG_{S}(T^*X;\mu_{\Delta_X}(h_{X*}\omega_X)) \ar[d]^{\wr}_{\Theta_{i_M}} \ar[r]^{\Theta_f} & \RG_{S^{\prime \prime}}(T^*Y;\mu_{\Delta_Y}(h_{Y*}\omega_Y)) \ar[d]^{\wr}_{\Theta_{i_N}} \\
\RG_{S\cap \Lb}(T^*M;\mu_{\Delta_M}(\delta_{M*}\omega_M)) \ar@{-}[d]^{\wr} \ar[r]^{\Theta_g} &\RG_{S^{\prime \prime} \cap \sLb}(T^*N;\mu_{\Delta_N}(\delta_{N*}\omega_N)) \ar@{-}[d]^{\wr} \\
\RG_{S\cap \Lb}(T^*M;\pi_M^{-1}\omega_M) \ar[r]^{g^*} & \RG_{S^{\prime \prime} \cap \sLb}(T^*N;\pi_N^{-1}\omega_N).}
\end{equation}
The bottom square in Diagram 8.b is obtained in this way. Moreover, by the construction of the morphism $\Theta_{i_M}$ (resp. $\Theta_{i_N}$) the image of $\id_F$ (resp. $\id_{f^{-1}F}$) by the left vertical arrows (resp. the right vertical arrows) in Diagram 8.b is $\overline{LC(F,\Phi_X)_M}$ (resp. $\overline{LC(f^{-1}F,\Phi_Y)_N}$). Hence the desired formula \eqref{eq:7-113} follows from the commutativity of Diagram 8.b. This completes the proof. \qed
\end{proof}

As a special case of this theorem, we obtain the following result which drops the condition \eqref{eq:6-42} of Theorem \ref{thm:6-4}.

\begin{cor}\label{cor:7-11}
In the situation as Theorem \ref{thm:6-4}, instead of assuming the condition \eqref{eq:6-42}, assume that the inclusion map $i_M \colon M \longhookrightarrow X$ of the fixed point manifold $M$ is non-characteristic for $F$. Then we have
\begin{equation}
LC(F,\Phi)_M =\sgn(\id-\phi^{\prime}) \cdot LC(F|_M,\Phi|_M)_M
\end{equation}
in $T^*M$. In particular, if $\supp(F) \cap M$ is compact, the local contribution $c(F,\Phi)_M$ of $(F,\Phi)$ from $M$ is expressed by the topological integral of the constructible function $\varphi(F|_M,\Phi|_M)$ on $M$:
\begin{equation}
c(F,\Phi)_M =\sgn(\id-\phi^{\prime}) \cdot \dint_M \varphi(F|_M,\Phi|_M).
\end{equation}
\end{cor}

\begin{rem}
Corollary \ref{cor:7-11} is not true if we do not assume that $i_M \colon M \longhookrightarrow X$ is non-characteristic for $F$. See e.g. \cite[Example 9.6.18]{K-S}.
\end{rem}

\end{document}